\documentclass[12pt]{article}
\usepackage{amsmath,amssymb,amsthm}

\usepackage{fullpage}

\title{Ramanujan's Most Singular Modulus}

\author{MARK B. VILLARINO \\[12pt]}

\newtheorem{theorem}{THEOREM}
\newtheorem{lemma}[theorem]{LEMMA}
\theoremstyle{definition}
\newtheorem{definition}{DEFINITION}
\newtheorem{example}{Example}

\numberwithin{equation}{section}

\newcommand{\dis}{\displaystyle} 
\newcommand{\hideqed}{\renewcommand{\qed}{}} 
\newcommand{\N}{\mathbb{N}}      
\newcommand{\phrase}[1]{\quad\mbox{#1}\quad} 
\newcommand{\Q}{\mathbb{Q}}      
\newcommand{\R}{\mathbb{R}}      
\newcommand{\sump}{\sideset{}{'}\sum} 
\newcommand{\Z}{\mathbb{Z}}      

\begin{document}

\maketitle

\begin{abstract}
We present an elementary self-contained detailed computation of Ramanujan's most famous singular
modulus, $k_{210}$, based on the Kronecker Limit Formula. 
\end{abstract}

\tableofcontents
 \section{The Singular Modulus \boldmath$\alpha$}
\subsection{Introduction}

In his second letter to G. H. \textsc{Hardy}, dated 27 February 1913,
the self-taught Indian genius \textsc{Srinivasa Ramanujan} announced the
following outrageous theorem (numbered (17) in the letter~\cite{rama}):

{\Large
\begin{theorem}[Ramanujan's Singular Modulus]
If
\begin{equation}
F(\alpha) := 
1 + \left(\frac{1}{2}\right)^2 \alpha
+ \left(\frac{1{\cdot}3}{2{\cdot}4}\right)^2 \alpha^2 +\cdots
\label{eq:Rama-thm-series}
\end{equation}
and
\begin{equation}
F(1 - \alpha) = \sqrt{210} \cdot F(\alpha),
\label{eq:Rama-thm-eqn}
\end{equation}
then
\begin{equation}
\fbox{$\dis \begin{aligned}
\alpha &= (4 - \sqrt{15})^4 (8 - 3\sqrt{7})^2 (2 - \sqrt{3})^2
(6 - \sqrt{35})^2
\\
&\quad \cdot (\sqrt{10} - 3)^4 (\sqrt{7} - \sqrt{6})^4
(\sqrt{15} - \sqrt{14})^2 (\sqrt{2} - 1)^2 .
\end{aligned}$}
\label{eq:Rama-thm-alpha}
\end{equation}
\end{theorem}
}

\normalsize

(We have written $\alpha$ in place of Ramanujan's $k$.)

Our paper develops an elementary self-contained presentation of the
history, theory and algorithms for computing $\alpha$ as well as the
details of its computation. Although
\eqref{eq:Rama-thm-series}--\eqref{eq:Rama-thm-alpha} ``...is one of
Ramanujan's most striking results,'' (Hardy~\cite{h}), there is no
published \textit{ab initio} account of it. \textsc{Watson} \cite{wa}
published the only computation of $\alpha$ in the literature, but it
presupposes a considerable background knowledge on the part of the
reader.

 There always has been and continues to be a large population of students and practicing mathematicians who are intrigued by the folklore surrounding \textsc{Ramanujan}'s origins, who are thrilled by the wondrous classical beauty of his mathematics, but who are put off by the necessity of mastering multiple branches of advanced mathematics in order to be able to read the proofs.
  
  We present for such readers  an expository development of one of his most famous formulae, a detailed elementary self-contained connected account which takes an inquiring mind with a modest background step-by-step to the full result.  Our exposition thus contains historical information, motivational commentary, and offers a ``cultural framework" as it were, so as to explain to the reader where the material ``fits in" within the greater subject of mathematics in general.  It is far more than just a computation.  It is a focused introduction to and an exposition of that specific area of mathematics which \textsc{Ramanujan}'s result so strikingly illuminates.
  
 Although our paper is expository, one hopes that the specialist can enjoy the point of view and some novelties of detail.  To begin with, our version of Ramanujan's ``very curious
algebraical lemma'' is new variant, and we use it for the major computations
in our paper (see \S 2.4.2).
    
Moreover, our proof of \textsc{Weber}'s fundamental formula for
$g_{2n}$ avoids Weber's use of genus theory by appealing to an
explicit composition formula for the quadratic form $AX^2 + BY^2$ (see
\S 3.7.1). This new approach simplifies the proof considerably and
brings it within the reach of our target reader.
   
This reader needs to understand elementary analysis and algebra at the
level of, say, three years at a university, as well as to be familiar
with number theory to the level of quadratic reciprocity and the
\textsc{Jacobi} symbol. All of this material is more than sufficiently
covered in \textsc{Niven} and \textsc{Zuckerman}~\cite{nz}.
   
Much of the underlying general theory has been published in
\textsc{David Cox}'s beautiful book \cite{cox} but he never deals with
Ramanujan's particular result. Finally, much of the material is in
scattered references in several languages.
 
We call special attention to Volume 5 of \textsc{Bruce Berndt}'s
edition of Ramanujan's Notebooks ~\cite{berndt}, which devotes more
than one hundred and fifty pages to the computation of singular
moduli, and which is the most exhaustive reference on the subject.
Unfortunately, although  Berndt presents \textsc{Chan}'s proof of the ``very curious algebraical lemma," he does not discuss the complex computation of $g_{210}$, the most difficult step, and thus is not a source for computing Ramanujan's specific result.  However, he does
 give some bibliographic references.

\subsection{Singular moduli and units in quadratic fields}

If we take Ramanujan's theorem at face value and try to
see if it makes sense (!) we can observe from the definition:
\begin{equation}
\fbox{$\dis
F(\alpha) :=
1 + \left(\frac{1}{2}\right)^2 \alpha
+ \left(\frac{1{\cdot}3}{2{\cdot}4}\right)^2 \alpha^2 +\cdots$}
\label{eq:Fa-series}
\end{equation}
that $F(\alpha)$ converges for $-1 \leq \alpha < 1$. Moreover, as
$\alpha$ increases from $0$ to $1$, $F(\alpha)$ \textit{increases}
from $1$ to infinity and $F(1 - \alpha)$ \textit{decreases} from
infinity to $1$, and therefore
$$
\frac{F(1 - \alpha)}{F(\alpha)}
\label{eq:Fa-ratio}
$$
\textit{decreases monotonically} from infinity to zero as $\alpha$
increases from $0$ to $1$. Thus \textit{there exists a unique real
number} $\alpha$ \textit{with} $0<\alpha<1$ \textit{such that}
\begin{equation}
\frac{F(1 - \alpha)}{F(\alpha)} = \sqrt{210}.
\end{equation}
The number
$$
k_{210} := \sqrt{\alpha}\equiv \sqrt{\alpha_{210}}
$$
is called the \textbf{singular modulus} (for 210). If the right hand
side of \eqref{eq:Fa-ratio} is $\sqrt{n}$, $n \in \R^+$, the same
argument shows that there exists a \textit{unique real number}
$\alpha_n$ with $0\leq\alpha_n<1$ which satisfies
$$
\frac{F(1 - \alpha_n)}{F(\alpha_n)} = \sqrt{n},
$$
and we write
$$
k_n := \sqrt{\alpha_n}.
$$
and we call $k_n$ the \textbf{singular modulus} (for $n$). We use
$\alpha$ and $k$ without the subscript quite frequently for ease of
exposition. Its meaning should be clear from the context. Of course,
what makes Ramanujan's result
\eqref{eq:Rama-thm-series}--\eqref{eq:Rama-thm-alpha} so astonishing
is the \textit{explicit} numerical value of $\alpha_{210}$ as a finite
product of the differences of quadratic surds. It is all the more
amazing and unexpected since there is no evident \textit{a priori}
reason to expect such a value based on the expansion
\eqref{eq:Fa-series}. \textit{Where does} $\alpha$ \textit{come from}?

Let's examine $\alpha$, itself. In the \textit{first} place, all the
quadratic surds,
$$
\sqrt{15}, \sqrt{7}, \sqrt{3}, \sqrt{35},
\sqrt{10}, \sqrt{6}, \sqrt{14}, \sqrt{2}
$$
which appear in $\alpha$ are the square roots of the \textit{divisors}
of $210$. This can hardly be an accident!

In the \textit{second} place, $\alpha$ is a product
$$
\fbox{$\dis \begin{aligned}
\alpha &:= (4 - \sqrt{15})^4 (8 - 3\sqrt{7})^2 (2 - \sqrt{3})^2
(6 - \sqrt{35})^2
\\
&\quad (\sqrt{10} - 3)^4 (\sqrt{7} - \sqrt{6})^4 
(\sqrt{15} - \sqrt{14})^2 (\sqrt{2} - 1)^2
\end{aligned}$}
$$
of terms of the form $T - U\sqrt{m}$ where
\begin{equation}
T^2 - m U^2 = 1
\label{eq:Pell-eqn}
\end{equation}
which is the famous \textsc{Pell}'s equation. Moreover the numbers
$u:=T - U\sqrt{m}$, where $T, U$ are rational integers, are
\textit{integers in the quadratic field} $Q(\sqrt{m})$. The
\textit{norm} of the quadratic integer $u$ is defined by the equation
$N(u):=T^2 - mU^2$. See~\cite{nz}. The equation \eqref{eq:Pell-eqn}
shows that $u$ has \textit{norm equal to one}, that is, $u$ a
\textit{unit} in $Q(\sqrt{m})$. Thus, each factor in $\alpha_{210}$ is
a unit in a subfield of $Q(\sqrt{m})$, and therefore $\alpha_{210}$
\textit{itself is a unit}! An accident?

Not only is this \textit{not} an accident, but rather is an instance
of a general theorem which \textsc{Chan} and \textsc{Huang} \cite{ch}
proved in 1997:

\begin{theorem}[Chan--Huang Unit Theorem]
If $n$ is of the form $4m + 2$, then $\alpha_n$ is a \textbf{unit}.
\qed
\end{theorem}

We will call this result the \textsc{Chan--Huang}--UNIT--THEOREM, or
``CHUT'' for short. In fact, the CHUT deals with $\alpha_n$ for
\textit{any}~$n$, but this suffices for now.

The fact that $\alpha_{210}$ (and therefore $k_{210}$) is a unit
suggests that the explanation of its value lies in the theory of
quadratic number fields. But the theory of quadratic fields is a
reformulation of the more classical theory of \textit{quadratic
forms}, and so we now turn to discussing $\alpha$ in terms of this
classical theory.

\subsection{Binary quadratic forms. The ``numeri idonei'' of Euler}

In 1776, \textsc{Euler}~\cite{e2} made the following observation:

\begin{theorem}
Let $n$ be an odd number not divisible by 3 nor 5 nor 7, but which is
properly represented by $x^2 + 210y^2$. If the equation
$$
x^2 + 210 y^2 = n
$$
has only ONE solution $(x,y)$ with $x,y\geq0$, then $n$ is a
\textbf{prime number}.
\qed
\end{theorem}

This property of ``$210$'' is shared by $64$ other known numbers, each
of which was called a \textit{``numerus idoneus''} (a ``convenient''
or ``useful'' number) by Euler. \textsc{Gauss} \cite{gauss}
conjectured that these $65$ numbers are \textit{all} of them, and
\textsc{Weinberger}~\cite{wei} proved that at most there may be one
more.

This common ``useful'' property, in our case, can be formulated as
follows in terms of the theory of binary quadratic forms:

\begin{theorem}
All the binary quadratic forms of determinant $210$ belong to $8$
genera and each \textbf{genus} contains \textbf{one class} of
equivalent forms.
\qed
\end{theorem}

We will define the technical terms later on.

\subsection{Elliptic modular functions and abelian extensions}
\label{sec:j-tau}

Our formulation of the ``useful'' property of $210$ can be stated in
terms of the theory of \textit{class fields} as follows. Let $\tau$ be
any complex number with positive imaginary part. Define
$$
j(\tau) := \frac{\left[1 + 240\dis\sum_{n=1}^{\infty}\sigma_3(n)
e^{2\pi i n \tau}\right]^3}
{{e^{2\pi i\tau}\dis\prod_{n=1}^{\infty}[1 - e^{2\pi in\tau}]^{24}}}
$$
where
$$
\sigma_{3}(n):=\sum_{d|n,\ d\geq1} d^3
$$
and ``$d|n$'' means ``$d$ divides $n$''.

In the theory of elliptic functions (\textsc{Borwein} and
\textsc{Borwein}~\cite{bb}) one proves
\begin{equation}
j = 64\frac{(1 - \alpha + \alpha^2)^3}{\alpha^2(1 - \alpha)^2}
\label{eq:j-of-alpha}
\end{equation}
where $\tau:=\sqrt{-n}$, and $\alpha\equiv\alpha_n$. (The reader
should take note of the two uses of the letter $n$, as a summation
index and also as the value whose square root appears in the
variable). Taking $n=210$ and using Ramanujan's value of $\alpha$
\eqref{eq:Rama-thm-alpha} as well as the equation
\eqref{eq:j-of-alpha}, some tedious algebra leads us to the value
\begin{center}
\framebox{\parbox{15cm}{$$
j(\sqrt{-210}) =
64\left\{4[\sqrt{3} + \sqrt{2}]^{12}[3\sqrt{14} + 5\sqrt{5}]^4
\left[\frac{\sqrt{7} + \sqrt{3}}{2}\right]^{12}
\left[\frac{(\sqrt{5} + 1)}{2}\right]^{12} + 1\right\}^3
$$
$$
\qquad \left([\sqrt{3} - \sqrt{2}]^{12}[3\sqrt{14} - 5\sqrt{5}]^4
\left[\frac{\sqrt{7} - \sqrt{3}}{2}\right]^{12}
\left[\frac{\sqrt{5} - 1}{2}\right]^{12}\right)
$$}}
\end{center}

We collect some facts about $j(\sqrt{-210})$. (See \textsc{Cox} \cite{cox}.)

\begin{enumerate}
    \item  $j(\sqrt{-210})$ is an \textit{algebraic integer}.
    \item  $j(\sqrt{-210})$ is a primitive element of the
\textit{ring class field} of $Q(\sqrt{-210})$.
    \item  $j(\sqrt{-210})$ satisfies a monic algebraic equation of
degree~$8$ with rational integer coefficients:
$$
\alpha^8 + a_1 \alpha^7 +\cdots+ a_8 = 0.
$$
The coefficients are \emph{enormous}! For example, the coefficient of
the seventh power of the variable is
$$
a_1 = - 3494487845306481075093315600749304691200
$$
while the constant term is
\begin{align*}
a_8 = \ 
& 75871693802713797386369191426742800771304395043277
\\
& 326055125100897851220991378671072700656000000000000. 
\end{align*}
(We thank \textsc{Anthony Varilly} for performing these computations
with the program PARI.)
    \item Combining this monic equation and the equation
\eqref{eq:j-of-alpha}, we conclude that $k_{210}$ satisfies a monic
equation of degree~$96$.
\end{enumerate}
This theory is called the theory of \textit{Complex Multiplication}
and today is in the forefront of modern mathematical research.

We will expand on these results later.

\subsection{Prospectus}

Let's stop to catch our breath. We've used several (as of yet)
undefined technical terms and notions. We've wandered over complex
function theory, elementary and algebraic number theory, elliptic and
modular functions, the genus theory of binary quadratic forms, and
finally, complex multiplication and class field theory (the theory of
normal field extensions with abelian Galois group)!

Did Ramanujan know any of this? Hardy \cite{h} suggests that he most
certainly did \textit{not} know the algebraic number theory, much less
class field theory.

Then \textit{how} did he arrive at his results? Nobody knows, really,
although Watson has made a clever suggestion which we will not
consider in this paper. (See \cite{wa} and \cite{chan}).

Our paper presents an elementary method of computing
$\alpha_{210}$.\footnote{The method suggested by \textsc{Watson} \cite{wa} and \textsc{Chan} \cite{chan} is superficially more elementary.  It entails ``guessing" the fields which contain the singular moduli and then using skillful ingenuity with radicals to obtain the final formulas.  Unfortunately, most readers have not mastered class-field theory, the theorems of \textsc{Chan} are not applicable to \textsc{Ramanujan}'s computation in any event, and \textsc{Watson}'s presentation presupposes a profound knowledge, as \textsc{Ramanujan} had, of singular moduli.  Indeed, we think that this is the method \textsc{Ramanujan} most likely used for the bulk of his computations of singular moduli, but it is \emph{not} remotely ``elementary" for most readers who don't know class field theory.} It is based on the theory of binary quadratic forms
and uses the so-called \textsc{Kronecker} \textit{Limit Formula}. We
will develop the details in full since they have never been published.
The fundamental sources are Weber's treatise \cite{w1} and his paper
\cite{w2}. In both references Weber simply \emph{states} the result of
his computation, but suppresses the computation itself. Ramanujan uses
the result of Weber's computation in order to compute $\alpha_{210}$,
although he obtained it independently of Weber. The fact is, nobody
knows how he did it.

Our plan is the following: we will first develop the theory of the
function $F(\alpha)$ and learn the \textit{two-step algorithm }which
Ramanujan used. Then we will develop the theory of quadratic forms and
the theory of Kronecker's Limit Formula so as to be able to carry out
the steps in Ramanujan's algorithm. Finally we will carry out the
details of the computation and crown our memoir with an elementary and
detailed demonstration of the general formula implied by the
computation which will explain a\textit{ priori} the miraculous
computations we carry out.

\section{Ramanujan's Function \boldmath$F(\alpha)$}
\subsection{Complete elliptic integrals}

We began our inquiry of Ramanujan's theorem
\eqref{eq:Rama-thm-series}--\eqref{eq:Rama-thm-alpha} by looking at
the value of $\alpha$. \textit{Now} we look at the function
$F(\alpha)$. In fact, this function is quite well known. Define
\begin{equation}
\fbox{$\dis
K(k) := \int_0^{ \frac{\pi}{2}}
\frac{d\theta}{\sqrt{1 - k^2\sin^2\theta}}$}
\label{cei} 
\end{equation}
Then $K(k)$, for $0\leq k<1$, is called the \textbf{complete
elliptic integral of the first kind}.  If we put
$$
\alpha := k^2,
$$
and then expand the denominator by the binomial theorem and integrate
term by term, \textit{we see that Ramanujan's} $F(\alpha)$
\textit{is} $K(k)$:
\begin{equation}
\fbox{$\dis F(\alpha) = K(k)$}
\label{eq:Fa-Kk}
\end{equation}
The number $k$ is called the \textbf{modulus} of $K$. The number $k'$,
defined by
$$
k'{}^2 := 1 - k^2,
$$
is called the \textbf{complementary modulus} and
\begin{equation}
K' := K(k')
\label{eq:compl-modl}
\end{equation}
is called the \textbf{complementary integral}. We also write
$$
K \equiv K(k).
$$

Mathematicians have studied $K$ and $K'$ for (at least) two and a half
centuries, and therefore they've studied Ramanujan's $F(\alpha)$ and
$F(1 - \alpha)$. See Euler \cite{e1} and
\textsc{Fagnano}~\cite{fagnano}.

\subsection{The modular equation}

In particular, in $1771$ \textsc{Landen}~\cite{landen} discovered the
following transformation law :
\begin{align}
K(k) &= \frac{1}{1+k} K\biggl(\frac{2\sqrt{k}}{1+k}\biggr),
\label{Landen}
\\[\jot]
K'(k) &= \frac{2}{1+k} K'\biggl(\frac{2\sqrt{k}}{1+k}\biggr).
\nonumber
\end{align}

Dividing $K'(k)$ by $K(k)$, we obtain
\begin{equation}
\frac{K'(k)}{K(k)}
= 2 \, \frac{\dis K'\biggl(\frac{2\sqrt{k}}{1+k}\biggr)}
{\dis K\biggl(\frac{2\sqrt{k}}{1+k}\biggr)},
\label{eq:dupl-law}
\end{equation}
which is a \textit{multiplication law} for the quotient
$\dfrac{K'}{K}(k) \equiv \dfrac{K'(k)}{K(k)}$. If we put
\begin{equation}
l := \frac{2\sqrt{k}}{1+k},
\label{eq:l-modl}
\end{equation}
the \textit{duplication law} \eqref{eq:dupl-law} takes the form
\begin{equation}
2 \, \frac{L'}{L} = \frac{K'}{K},
\label{eq:dupl-law-bis}
\end{equation}
where $L'$ and $L$ are the complete elliptic integrals with respective
moduli $l'$ and~$l$. The equation \eqref{eq:l-modl} can also be
written
\begin{equation}
l^2(1 + k)^2 = 4k.
\label{eq:modl-eqn-two}
\end{equation}
Equations \eqref{eq:l-modl} and \eqref{eq:modl-eqn-two} are both
called the \textbf{modular equation of degree 2}. The duplication law
\eqref{eq:dupl-law-bis} is one of an infinity of multiplication laws
and each such law has a modular equation associated with it. We now
make the following general definition.

\begin{definition}
The \textbf{\boldmath modular equation of degree~$n$} is the algebraic
equation
$$
\Omega_n(k,l) = 0,
$$
relating the moduli $k$ and $l$ for which the multiplication law
\begin{equation}
n \frac{L'}{L} = \frac{K'}{K}
\label{eq:multn-law}
\end{equation}
holds. Here $K, K', L, L'$ are the complete elliptic integrals with
moduli $k$ and~$l$, and $n > 0$ is a rational number.
\end{definition}

Jacobi~\cite{jacobi} proved \textit{that for every integer} $n>0$
\textit{there is a modular equation} $\Omega_n(k,l)=0$. We just proved
Jacobi's result, with Landen's help (!) for $n=2$. We list a few more
examples of modular equations (\textsc{Greenhill} \cite{gh}):

\begin{center}
\begin{tabular}{|c|c|}
\hline
\rule[-3mm]{0mm}{8mm}\textbf{Degree}& \textbf{Modular Equation}\\
\hline \hline
\rule[-3mm]{0mm}{8mm}3 & $\sqrt{kl} + \sqrt{k'l'}=1$ \\
\hline
\rule[-3mm]{0mm}{8mm}5 & $kl + k'l' + \sqrt[3]{32klk'l'}=1$ \\
\hline
\rule[-3mm]{0mm}{8mm}7 & $\sqrt[4]{kl} + \sqrt[4]{k'l'}=1$ \\
\hline
\end{tabular}
\end{center}

\subsection{Ramanujan's theorem for \boldmath$\alpha_2$, $\alpha_3$,
and $\alpha_7$}
\label{ssc:small-alpha}

\textit{Why} are we interested in the modular equation? Equation
\eqref{eq:multn-law} is
$$
n \frac{K(l')}{K(l)} = \frac{K(k')}{K(k)}
$$
Now comes the ``brilliancy'', as they say in chess. \textit{Choose}
$k$ \textit{such that}
$$
\fbox{$l = k',  \quad  l' = k.$}
$$
Then \eqref{eq:multn-law} becomes
$$
n \frac{K}{K'} = \frac{K'}{K},
$$
or
$$
\Bigl( \frac{K'}{K} \Bigr)^2 = n,
$$
that is,
\begin{equation}
\frac{K'}{K}(k) = \sqrt{n}.
\label{eq:Rama-eqn-k}
\end{equation}
or, by \eqref{eq:Fa-Kk}, \eqref{eq:compl-modl} and \eqref{eq:Rama-eqn-k},
\begin{equation}
\fbox{$\dis \frac{F(1-\alpha)}{F(\alpha)} = \sqrt{n}$,}
\label{eq:Rama-eqn}
\end{equation}
and \eqref{eq:Rama-eqn} \textit{is precisely Ramanujan's equation}
\eqref{eq:Rama-thm-eqn} ---indeed it is for \textit{any}~$n$, and not
just $n = 210\,$!

So we now understand the origin of Ramanujan's equation in general:
namely
\begin{center}
\fbox{\parbox{8cm}{\itshape The modular equation of degree~$n$ takes
the form of Ramanujan's equation \eqref{eq:Rama-eqn} when in it we
choose $k = l'$ and $k' = l$.}}
\end{center}

Let us compute a few examples.

\begin{example} 
\fbox{$n = 2$}\quad
Suppose $\dis \frac{F(1 - \alpha)}{F(\alpha)} = \sqrt{2}$. Then the 
modular equation of degree $2$, namely
$$ 
l^2(1 + k)^2 = 4k,
$$
under the substitution $l = k'$, $l' = k$, becomes
$$
(1 - k^2)(1 + k)^2 = 4k,
$$
or
$$
(k^2 + 1) (k^2 + 2k - 1) = 0,
$$
or
$$
k = -1 \pm \sqrt{2}.
$$
But $0 \leq k_2 \leq 1$, and therefore
$$
k_2 = \sqrt{2} - 1.
$$
\end{example}

\begin{theorem}
If $F(\alpha) = \sqrt{2} \cdot F(1 - \alpha)$, then
\boldmath
$$
\alpha = (\sqrt{2} - 1)^2.
\eqno \qed
$$
\end{theorem}

Thus we have found the singular modulus $k_2$.

\begin{example} 
\fbox{$n = 3$}\quad
Our table of modular equations gives us the following equation of
degree~$3$:
$$
\sqrt{kl} + \sqrt{k'l'} = 1.
$$
Putting $l = k'$, $l' = k$, we obtain 
$$
2\sqrt{kk'} = 1,
$$
or
$$
k^2 (1 - k^2) = \frac{1}{16},
$$
that is,
$$
k = \sqrt{\frac{2 \pm \sqrt{3}}{4}}.
$$
To properly choose the sign, we recall that
$$
\alpha = k^2 = \frac{2 \pm \sqrt{3}}{4}
$$
(which makes $1 - \alpha = (2 \mp \sqrt{3})/4$), is the
\textit{unique} solution in $[0,1]$ to the equation
\begin{align*}
\sqrt{3} &= \dfrac{F(1 - \alpha)}{F(\alpha)}
\\
&= 
\dfrac{1 + \Bigl(\dfrac{1}{2}\Bigr)^2
  \biggl[\dfrac{2 \mp \sqrt{3}}{4}\biggr]
  + \Bigl( \dfrac{1 \cdot 3}{2 \cdot 4} \Bigr)^2 
  \biggl[\dfrac{2 \mp \sqrt{3}}{4}\biggr]^2 + \cdots}
{1 + \Bigl(\dfrac{1}{2}\Bigr)^2
  \biggl[\dfrac{2 \pm \sqrt{3}}{4}\biggr]
  + \Bigl( \dfrac{1 \cdot 3}{2 \cdot 4} \Bigr)^2 
  \biggl[\dfrac{2 \pm \sqrt{3}}{4}\biggr]^2 + \cdots}
\end{align*}
and comparing both sides numerically, we see that we must choose the
\textit{negative} sign. Therefore
$$
k_3 = \sqrt{\frac{2 - \sqrt{3}}{4}}.
$$
and our Ramanujan type theorem is the following.
\end{example}

\begin{theorem}
If $F(\alpha) = \sqrt{3} \cdot F(1 - \alpha)$, then
\boldmath
$$
\alpha = \frac{2 - \sqrt{3}}{4}.
\eqno \qed
$$
\end{theorem}

We have therefore found the singular modulus $k_{3}$.
We note that the choice
$$
\alpha = \frac{2 + \sqrt{3}}{4}
$$
is the root of the \textit{complementary} equation
$$
F(1 - \alpha) = \frac{1}{\sqrt{3}} \, F(\alpha).
$$

\begin{example} 
\fbox{$n = 7$}\quad
The modular equation of degree~$7$ (by our table) is:
$$
\sqrt[4]{kl} + \sqrt[4]{k'l'} = 1.
$$
Putting $l = k'$, $l' = k$, and solving for $k_7$, we obtain
$$
k_7 = \sqrt{\frac{8 - \sqrt{7}}{16}}.
$$
(We leave the details to the reader.)
\end{example}

\begin{theorem}
If $F(\alpha) = \sqrt{7} \cdot F(1 - \alpha)$, then
\boldmath
$$
\alpha = \frac{8 - 3\sqrt{7}}{16}.
\eqno \qed
$$
\end{theorem}

We have therefore found the singular modulus $k_{7}$.


\subsection{Ramanujan's two-step algorithm}

The computation of $\alpha_2$, $\alpha_3$, and $\alpha_7$ we just
carried out would suggest a procedure to calculate Ramanujan's
$\alpha_{210}$. Namely, \textit{find the modular equation of degree
$210$ and then transform it into an algebraic equation for $k_{210}$
by means of the substitution $k=l'$ and $k'=l$}.

Although it is \textit{theoretically} possible to construct
$\Omega_{210}(k,l)$, the labor is formidable beyond belief and the
coefficients are so huge and unwieldly as to make it impossible in
practice. We just have to recall the coefficients of the equation for
$j(\sqrt{-210})$ that we saw in Section~\ref{sec:j-tau}. Ramanujan
unquestionably did \textit{not} follow this procedure.

Then what \textit{did} he do?

Fortunately, he indicated the general procedure; but unfortunately, he
left no clue about how to carry out the hardest part.

\subsubsection{Ramanujan's invariant \boldmath$g_n$}

One of Jacobi's great discoveries is the \emph{infinite product
representation} of the moduli $k$ and $k'$ (see \cite{bb}). Let
$$
q := e^{-\pi K'/K}.
$$
Then Jacobi's marvelous formulas are:
\begin{equation}\begin{split}
k(q) &= 4\sqrt{q} 
\biggl[ \frac{1 + q^{2}}{1 + q}\ \frac{1 + q^{4}}{1 + q^{3}}\ \frac{1 + q^{6}}{1 + q^{5}}\ \frac{1 + q^{8}}{1 + q^{7}}\ \cdots \biggr]^4\\
&\\
k'(q) &= 
\biggl[ \frac{1 - q}{1 + q}\ \frac{1 - q^{3}}{1 + q^{3}}\ \frac{1 - q^{5}}{1 + q^{5}}\ \frac{1 - q^{7}}{1 + q^{7}}\ \cdots \biggr]^4.\end{split}
\label{eq:Jacobi-prods}
\end{equation}
Needless to say, Ramanujan was intimately familiar with both of these
formulas. Indeed, it is highly probable that he (re)discovered them!
If we assume $\dfrac{K'}{K}(k) = \sqrt{n}$, then we follow
Ramanujan~\cite{rama} and \textit{we define} the function $g_n$ by the
equation:
{\large \begin{equation}
\fbox{$\dis
\mbox{\boldmath $g_n$} := 2^{-\frac{1}{4}} e^{\frac{\pi\sqrt{n}}{24}}
(1 - e^{-\pi\sqrt{n}}) (1 - e^{-3\pi\sqrt{n}}) (1 - e^{-5\pi\sqrt{n}})
\cdots
$}
\label{eq:Rama-prod}
\end{equation}}
where we have used boldface type for emphasis. Then
\eqref{eq:Jacobi-prods} and \eqref{eq:Rama-prod} and some algebra
permit us to conclude that
\begin{equation}
\frac{k'^2}{k} = 2g_n^{12},
\label{eq:kkk-gn}
\end{equation}
and solving \eqref{eq:kkk-gn} for the singular modulus $k\equiv k_n$,
we obtain the formula
\begin{equation}
k_n = g_n^6 \left(\sqrt{g_n^{12} + \frac{1}{g_n^{12}}} - g_n^6\right).
\label{eq:gggg-kn}
\end{equation}
Watson~\cite{wa} was the first person to point out \textbf{Ramanujan's
TWO-STEP ALGORITHM:}\footnote{We do not claim that \textsc{Ramanujan} was the first person to use this algorithm; we only claim that this was the one he used.}
\begin{center}
\fbox{\fbox{\parbox{5cm}{\begin{itemize}
\item Step 1: compute $g_n\,$,
\item Step 2: compute $k_n\,$.
\end{itemize}}}}
\end{center}

Because of \eqref{eq:gggg-kn}, \textit{all the difficulty is
concentrated in the computation of~$g_n$}.

\subsubsection{``a very curious algebraical lemma''}

If we recall our earlier computation of $k_2$, $k_3$, and $k_7$ in
Section~\ref{ssc:small-alpha}, we note that they were all
representable as products of \textit{quadratic units}, i.e., of
quantities of the form $\sqrt{A} - \sqrt{A-1}$. Similar computations
give us the values
\begin{align}
k_6 &= (2 - \sqrt{3})(\sqrt{3} - \sqrt{2}),
\nonumber \\
k_{10} &= (\sqrt{10} - 3)(3 - 2\sqrt{2}),
\nonumber \\
k_{30} &=
(5 - 2\sqrt{6})(4 - \sqrt{15})(\sqrt{6} - \sqrt{5})(2 - \sqrt{3}),
\label{eq:k-thirty}
\end{align}
which are again all products of quadratic units, as is Ramanujan's
value of $k_{210}$, a fact noted earlier.

Unfortunately, the formula \eqref{eq:gggg-kn}:
$$ 
k_n = g_n^6 \left(\sqrt{g_n^{12} + \frac{1}{g_n^{12}}} - g_n^6\right)
$$
does not easily transform itself into a product of quadratic units.
Perhaps the reader would like to try to transform \eqref{eq:gggg-kn}
into \eqref{eq:k-thirty} by using the value
$$
g_{30}^6 = (3 + \sqrt{10}) (2 + \sqrt{5}).
$$

Observe that $g_{30}$ is a \textit{product} of units.  The same
is true of all $g_n$ used by Ramanujan. Writing \eqref{eq:kkk-gn} in
the form
\begin{equation}
\frac{1}{k} - k = 2 g_n^{12},
\label{eq:kk-gntwelve}
\end{equation}
we see that Ramanujan \textit{had to solve the quadratic equation
\eqref{eq:kk-gntwelve}, in which $g_n$ is a \emph{product} of units,
and then express the \emph{root} $k_n$ again as a \emph{product} of
units}.

We now quote Hardy~\cite{h}:
\begin{center}
\fbox{\parbox{9cm}{``This result Ramanujan [achieved], granted the
value of [$g_n$], by a very curious algebraical lemma.''}}
\end{center}

For brevity, we refer to this lemma as the ``VCAL.'' We now reproduce
Watson's statement~\cite{wa} of the VCAL, which he, in turn, copied
from Ramanujan's notebook (with a few changes in notation):

\begin{theorem}[VCAL]
Suppose that
\begin{enumerate}
\item $uv := g_n^6$,
\item $2U := u^2 + \dfrac{1}{u^2};$\quad $2V := v^2 + \dfrac{1}{v^2}$,
\item $W := \sqrt{U^2 + V^2 - 1}$, \ and
\item $2S := U + V + W + 1$;
\end{enumerate}
then
\begin{align}
\alpha &= (\sqrt{S} - \sqrt{S-1})^2 (\sqrt{S-U} - \sqrt{S-U-1})^2
\nonumber \\
&\quad \cdot (\sqrt{S-V} - \sqrt{S-V-1})^2
(\sqrt{S-W} - \sqrt{S-W-1})^2.
\tag*{\qed}
\end{align}
\end{theorem}

The VCAL shows that $\alpha_n$ (and therefore $k_n$) is a product of
units in certain algebraic fields. Watson proved the VCAL by
verification without any indication of its origin. A more natural
proof may be found in Berndt~\cite{berndt}.

We have found the following \textit{alternate} version of the VCAL,
which we believe to be new, though it can be obtained from the work of \textsc{Chan} \cite{chan}, the ``AVCAL'', to be useful.

\begin{theorem}[AVCAL]
If
\begin{align}
\sqrt{\alpha} &:= \sqrt{ab} + \sqrt{(a+1)(b-1)},
\label{eq:avcal-alpha}
\\
\sqrt{\beta} &:= \sqrt{cd} + \sqrt{(c-1)(d-1)},
\nonumber
\end{align}
then
\begin{align*}
x_1 &:= (\sqrt{a+1} - \sqrt{a}) (\sqrt{b} - \sqrt{b-1}) 
(\sqrt{c} - \sqrt{c-1}) (\sqrt{d} - \sqrt{d-1}),
\\
x_2 &:= - (\sqrt{a+1} + \sqrt{a}) (\sqrt{b} + \sqrt{b-1})
(\sqrt{c} + \sqrt{c-1}) (\sqrt{d} + \sqrt{d-1})
\end{align*}
are the roots of
$$
\frac{1}{x} - x
= 2 \bigl[ \sqrt{\alpha\beta} + \sqrt{(\alpha + 1)(\beta - 1)} \bigr].
\eqno \qed
$$
\end{theorem}

The proof is quite simple. We offer the following two examples.

\subsubsection{The computation of \boldmath$\alpha_{30}$}
   
{}From Ramanujan's paper ``Modular equations and approximations
to~$\pi$'' \cite[p.~26]{rama}, we find that
\begin{align*}
g_{30}^{12} &= (2 + \sqrt{5})^2 (3 + \sqrt{10})^2
\\
&= 171 + 54\sqrt{10} + 76\sqrt{5} + 120\sqrt{2}.
\end{align*}
Take
\begin{align*}
\sqrt{\alpha\beta} &:= 171 + 54\sqrt{10},
\\
\sqrt{(\alpha + 1)(\beta - 1)} &:= 76\sqrt{5} + 120\sqrt{2}.
\end{align*}
Then
\begin{align*}
\alpha &= 759 + 240\sqrt{10} = (8\sqrt{6} + 5\sqrt{15})^2
\\
\implies \sqrt{\alpha} &= 8\sqrt{6} + 5\sqrt{15},
\end{align*}
and
\begin{align*}
\beta &= 39 + 12\sqrt{10} = (2\sqrt{6} + \sqrt{15})^2
\\
\implies \sqrt{\beta} &= 2\sqrt{6} + \sqrt{15}.
\end{align*}
Our value of $\alpha$ and \eqref{eq:avcal-alpha} suggest that we take
\begin{align*}
\sqrt{ab} &:= 8\sqrt{6}
\\
\implies ab &= (8\sqrt{6})^2 = 384,
\end{align*}
and
\begin{align*}
\sqrt{(a + 1)(b - 1)} &:= (5\sqrt{15})
\\
\implies (a + 1)(b - 1) &= (5\sqrt{15})^2 = 375.
\end{align*}
Solving for $a$ and $b$, we obtain
$$
a = 24,  \quad  b = 16.
$$
Similarly,
$$
cd = 24,  \quad  (c - 1)(d - 1) = 15,
$$
and therefore
$$
c = 6,  \quad  d = 4.
$$
Since $0 \leq \alpha < 1$, we conclude from the AVCAL that
$$ 
\fbox{$
k_{30} =
(5 - 2\sqrt{6}) (4 - \sqrt{15}) (\sqrt{6} - \sqrt{5}) (2 - \sqrt{3})
$}
$$
as given in Berndt's paper ``Ramanujan's singular moduli''
\cite[Thm.~2.1]{bcz}.

\subsubsection{The computation of \boldmath$k_{210}$}
   
This is the computation that inspired this paper. From Weber's tables
\cite{w1,w2}, or Watson's paper \cite{wa}, we learn that
\begin{align*}
g_{210}^{12}
&= \left( \sqrt{\sqrt{3} + \sqrt{2}}
\cdot \sqrt[6]{3\sqrt{14} + 5\sqrt{5}}
\cdot \sqrt{\frac{\sqrt{7} + \sqrt{3}}{2}}
\cdot \sqrt{\frac{\sqrt{5} + 1}{2}} \right)^{12}
\\
&= 120134025 + 53725540\sqrt{5} + 26215380\sqrt{21} + 11723880\sqrt{105}
\\
&\qquad + 49044510\sqrt{6} + 32107152\sqrt{14} + 21933360\sqrt{30}
+ 14358762\sqrt{70}.
\end{align*}
Now take
\begin{align*}
\sqrt{\alpha\beta} &:=
120134025 + 53725540\sqrt{5} + 26215380\sqrt{21} + 11723880\sqrt{105},
\\
\sqrt{(\alpha + 1)(\beta - 1)} &:=
49044510\sqrt{6} + 32107152\sqrt{14} + 21933360\sqrt{30}
+ 14358762\sqrt{70}.
\end{align*}
Solving for $\alpha$ and $\beta$, we obtain
\begin{align*}
\alpha &=
120621959 + 53943744\sqrt{5} + 26321856\sqrt{21} + 11676456\sqrt{105}
\\
&= (3168\sqrt{3} + 2076\sqrt{7} + 1419\sqrt{15} + 928\sqrt{35})^2
\\
\implies \sqrt{\alpha} &=
3168\sqrt{3} + 2076\sqrt{7} + 1419\sqrt{15} + 928\sqrt{35},
\\[2\jot]
\beta &=
119648071 + 53508216\sqrt{5} + 26109336\sqrt{21} + 11676456\sqrt{105}
\\
&= (3156\sqrt{3} + 2068\sqrt{7} + 1413\sqrt{15} + 924\sqrt{35})^2\\
\implies \sqrt{\beta} &=
3156\sqrt{3} + 2068\sqrt{7} + 1413\sqrt{15} + 924\sqrt{35}.
\end{align*}
Now since
\begin{align*}
\alpha &= \sqrt{ab} + \sqrt{(a + 1)(b - 1)}
\\
&= 3168\sqrt{3} + 2076\sqrt{7} + 1419\sqrt{15} + 928\sqrt{35},
\end{align*}
we take
\begin{align*}
\sqrt{ab} &:= 2076\sqrt{7} + 1419\sqrt{15},
\\
\sqrt{(a + 1)(b - 1)} &:= 3168\sqrt{3} + 928\sqrt{35}.
\end{align*}
Solving for $a$ and $b$, we obtain
$$
a = 121983 + 11904\sqrt{105},  \quad  b = 249 + 24\sqrt{105}.
$$
Similarly, in
\begin{align*}
\sqrt{\beta} &= \sqrt{cd} + \sqrt{(c - 1)(d - 1)}
\\
&= 3156\sqrt{3} + 2068\sqrt{7} + 1413\sqrt{15} + 924\sqrt{35},
\end{align*}
we take
\begin{align*}
\sqrt{cd} &:= 2068\sqrt{7} + 1413\sqrt{15},
\\
\sqrt{(c - 1)(d - 1)} &:= 3156\sqrt{3} + 924\sqrt{35}.
\end{align*}
Solving for $c$ and $d$, we obtain
$$
c = 121489 + 11856\sqrt{105},  \quad  d = 247 + 24\sqrt{105}.
$$
Therefore our AVCAL gives us the following formula for the singular
modulus $k_{210}$:
\begin{align*}
k_{210} &= \left[ \sqrt{121984 + 11904\sqrt{105}}
- \sqrt{121983 + 11904\sqrt{105}}\ \right]
\\[\jot]
&\qquad \cdot \left[ \sqrt{249 + 24\sqrt{105}}
- \sqrt{248 + 24\sqrt{105}}\ \right]
\\[\jot]
&\qquad \cdot \left[\sqrt{121489 + 11856\sqrt{105}}
- \sqrt{121488 + 11856\sqrt{105}}\ \right]
\\
&\qquad \cdot \left[\sqrt{247 + 24\sqrt{105}}
- \sqrt{246 + 24\sqrt{105}}\ \right].
\end{align*}

Now we follow Watson's computation \cite{wa}. The following identities
hold:
\begin{align*}
121984 + 11904\sqrt{105} &= (248 + 24\sqrt{105})^2 ,
\\
121983 + 11904\sqrt{105} &= (93\sqrt{7} + 64\sqrt{15})^2 ,
\\
249 + 24\sqrt{105} &= (12 + \sqrt{105})^2 ,
\\
248 + 24\sqrt{105} &= (6\sqrt{3} + 2\sqrt{35})^2 ,
\\
121489 + 11856\sqrt{105} &= (247 + 24\sqrt{105})^2 ,
\\
121488 + 11856\sqrt{105} &= (78\sqrt{10} + 38\sqrt{42})^2 ,
\\
247 + 24\sqrt{105} &= (4\sqrt{4} + 3\sqrt{15})^2 ,
\\
246 + 24\sqrt{105} &= (3\sqrt{14} + 2\sqrt{30})^2 .
\end{align*}
So also do the identities
\begin{align*}
248 + 24\sqrt{105} - (93\sqrt{7} + 64\sqrt{15})
&= (31 - 8\sqrt{15})(8 - 3\sqrt{7})
\\
&= (4 - \sqrt{15})^2(8 - 3\sqrt{7}) ,
\\
12 + \sqrt{105} - (6\sqrt{3} + 2\sqrt{35})
&= (6 - \sqrt{35})(2 - \sqrt{3}) ,
\\
247 + 24\sqrt{105} - (78\sqrt{10} + 38\sqrt{42})
&= (13 - 2\sqrt{42})(19 - 6\sqrt{10})
\\
&= (\sqrt{7} - \sqrt{6})^2(\sqrt{10} - 3)^2 ,
\\
4\sqrt{4} + 3\sqrt{15} - (3\sqrt{14} + 2\sqrt{30})
&= (3 - 2\sqrt{2})(\sqrt{15} - \sqrt{14})
\\
&= (\sqrt{2} - 1)^2(\sqrt{15} - \sqrt{14}) .
\end{align*}
Taken together, they transform our formula for $k_{210}$ into
{\Large $$
\fbox{$ \begin{aligned}
k_{210} = \ {}
& (4 - \sqrt{15})^2 (8 - 3\sqrt{7}) (6 - \sqrt{35}) (2 - \sqrt{3})
\\
& (\sqrt{7} - \sqrt{6})^2 (\sqrt{10} - 3)^2 (\sqrt{2} - 1)^2
(\sqrt{15} - \sqrt{14})
\end{aligned}
$}
$$}
\emph{which \textbf{is} the value that Ramanujan announced!}

\section{Kronecker's Limit Formula and the Computation of
\boldmath$g_{210}$}

Let's review our progress up to now. Our task is to compute $k_{210}$
by use of Ramanujan's two-step algorithm. We have just computed it,
\textit{granted} the value of $g_{210}$. Thus we still have to compute
the value of $g_{210}$.

There are at least four methods available for computing $g_n$:
\begin{enumerate}
    \item The Kronecker Limit Formula,
    \item Watson's empirical process,
    \item Modular equations,
    \item Class field theory.
\end{enumerate}

In this paper we will deal with Kronecker's limit formula. We begin
with an introductory sketch of the theory of quadratic forms, using the form $x^{2}+\mathbf{210}y^{2}$, fundamental to the computation of \textsc{Ramanujan}'s result, as the underlying example, then we
state the formula itself, and its specialization to the computation of
Ramanujan's function $g_n$, and finally we carry out the detailed
computation of the numerical value of $g_{210}$. Since the formula is
a statement about quadratic forms, we must first introduce the basic
concepts and terminology of this beautiful classical theory.

\subsection{Binary quadratic forms}

\begin{definition}
A \textbf{binary quadratic form} $F$ (or just a \textbf{form}) is a
polynonmial in two variables $X$ and $Y$, given by
$$
\fbox{$F \equiv F(X,Y) := aX^2 + bXY + cY^2$}
$$
with real coefficients $a$, $b$, and~$c$.
\begin{enumerate}
\item
$F$ is \textbf{rational} iff $a$, $b$, and $c$ are in~$\Q$,
\item
$F$ is \textbf{integral} iff $a$, $b$, and $c$ are in~$\N$,
\item
the \textbf{discriminant} of $F$, $\Delta$, is defined by
$$
\fbox{$\Delta \equiv \Delta(F) := b^2 - 4ac$,}
$$
\item
$F$ is \textbf{positive definite} iff
$$
\Delta < 0,  \quad  a > 0,  \phrase{and}  c > 0.
$$
\end{enumerate}
\end{definition}

\begin{definition}
An integral form $F(X,Y) := aX^2 + bXY + cY^2$
\textbf{represents} an integer~$n$ iff the equation
$$
aX^2 + bXY + cY^2 = n
$$
has a solution in integers $X$ and~$Y$.
\end{definition}

\begin{example}
\label{xmp:one}
The polynomial
$$
F(X,Y) := X^2 + 210 Y^2
$$
is a positive definite integral binary quadratic form of discriminant
$\Delta = -840$, and any number represented by $F(X,Y)$ is also
represented by the form
$$
G(X,Y) := 5266 X^2 - 8424 XY + 3369 Y^2,
$$
which again has discriminant $\Delta(G) = -840$. However, even though
the form
$$
H(X,Y) := 14 X^2 + 15 Y^2
$$
also is positive definite, integral, binary, quadratic and has
discriminant $\Delta(H) = -840$, it turns out that $H(X,Y)$ represents
\textit{no} number represented by $F$ and $G$, and \textit{vice
versa}.
\end{example}

The explanation of this example is to be found in the notion of
\textit{equivalence}.

\begin{definition}
$GL_2(\Z)$ is the multiplicative group of $2 \times 2$ matrices
$$
g := \begin{pmatrix} r & s \\ t & u \end{pmatrix}
$$
such that $r$, $s$, $t$ and $u$ are integers which satisfy
$ru - st = \pm 1$.

$SL_2(\Z)$ is the subgroup of matrices in $GL_2(\Z)$ with
$ru - st= +1$.
\end{definition}

\begin{example}
The matrix $g = \begin{pmatrix} 4 & -5 \\ -3 & 4 \end{pmatrix}$ lies
in $SL_2(\Z)$.
\end{example}

\begin{definition}
Given a form $F(X,Y) := aX^2 + bXY + cY^2$ and a matrix
$$
g = \begin{pmatrix} r & s \\ t & u \end{pmatrix} \in SL_2(\Z),
$$
we define the form $gF$ by the formula
$$
gF := a(rX + tY)^2 + b(rX + tY)(sX + uY) + c(sX + uY)^2.
$$
That is, $gF$ is got from $F$ by making the \textit{substitution}
$$
X \mapsto rX + tY,  \quad  Y \mapsto sX + uY.
$$
\end{definition}

\begin{example}
\label{xmp:three}
The transform of
$$
F(X,Y) := X^2 + 210 Y^2  \phrase{under}
g = \begin{pmatrix} 4 & -5 \\ -3 & 4 \end{pmatrix}
$$
is
\begin{align*}
gF &:= (4X - 3Y)^2 + 210(-5X + 4Y)^2
\\
&= 5266 X^2 - 8424 XY + 3369 Y^2
\\
&\equiv G(X,Y),
\end{align*}
and $G(X,Y)$ is the form introduced in Example~\ref{xmp:one}.
\end{example}

\begin{definition}
Two forms $F$ and $F'$ are \textbf{equivalent} iff there
exists a $g\in GL_2(\Z)$ such that $F' = gF$.

If there exists $g \in SL_2(\Z)$ such that $F' = gF$, we say that $F$
and $F'$ are \textbf{properly equivalent}.
\end{definition}

\begin{example}
The forms $F(X,Y)$ and $G(X,Y)$ of Example~\ref{xmp:one} and
Example~\ref{xmp:three} are properly equivalent.
\end{example}

We assume as known the following facts about quadratic forms. The
details may be found in \textsc{Flath}~\cite{flath}.
\begin{enumerate}
\item
Integral forms that are equivalent represent precisely
\textit{the same integers}.
\item
There is only a \textit{finite number, $h \equiv h(\Delta)$, of proper
equivalence classes of positive definite integral forms of a given
discriminant~$\Delta$}. This number $h$ is called the \textbf{class
number}.
\item
A positive definite form $F(X,Y) := aX^2 + bXY + cY^2$ is
called \textbf{reduced} iff
\begin{itemize}
\item
$|b| \leq a \leq c$,
\item
$|b| = a \implies b = a$, and
\item
$a = c \implies b \geq 0$.
\end{itemize}
\item
\textit{There is a \emph{unique} reduced form in every proper
equivalence class of positive definite integral forms.}
\end{enumerate}

\begin{example}
The form $x^{2}+\mathbf{210}y^{2}$ has discriminant $\Delta = -840$. Let us compute \textbf{\emph{all}} of the reduced forms of discriminant 
$\Delta = -840$.

The following Lemma is useful in numerical work.

\begin{lemma}
If $F$ is reduced and $\Delta < 0$, then
$$
a \leq \sqrt{\frac{|\Delta|}{3}}.
$$
\end{lemma} 

\begin{proof}
\begin{align}
4a^2 &\leq 4ac = b^2 - \Delta \leq a^2 - \Delta
\nonumber \\
\implies 3a^2 &\leq - \Delta = |\Delta|
\nonumber
\\
\implies a &\leq \sqrt{\frac{|\Delta|}{3}}.
\tag*{\qed}
\end{align}
\hideqed
\end{proof}

In our case, $a \leq \sqrt{840/3} = 16.7\dots$ and the candidates for
$a$ are $a = 1, 2, 3,\dots, 16$.

\vspace{6pt}
\noindent
Case $a = 1$:\quad
Then $|b| = 0,1$ since $|b| \leq a$. If $|b| = 0$, then
$$
c = \frac{b^2 + 840}{4a} = 210,
$$
and therefore we obtain the reduced form
$$
F(X,Y) = X^2 + 210 Y^2.
$$
If $|b| = 1$, then $c = (1 + 840)/4$ is not an integer, and we have
exhausted all possibilities with $a = 1$.

\vspace{6pt}
\noindent
Case $a = 2$:\quad
Then $|b| = 0, 1, 2$. Now $|b| = 0 \implies c = 840/8 = 105$, and
therefore
$$
F(X,Y) = 2 X^2 + 105 Y^2.
$$
On the other hand, $|b| = 1 \implies c = 841/8$, not an integer, and
$|b| = 2 \implies c = 844/8$, not an integer either.

The remaining values of $a$ are treated in the same way and we simply
tabulate the final results.
\end{example}

\begin{theorem}
The reduced forms with $\Delta = -840$ are:
\begin{center}
\framebox{\parbox{4cm}{\begin{enumerate}
    \item $\ \ X^2 + 210Y^2$\\
    \item $2X^2 + 105Y^2$\\
    \item $\ 3X^2 + 70Y^2$\\
    \item $\ 5X^2 + 42Y^2$\\
    \item $\ 6X^2 + 35Y^2$\\
    \item $\ 7X^2 + 30Y^2$\\
    \item $10X^2 + 21Y^2$\\
    \item $14X^2 + 15Y^2$
\end{enumerate}}}
\end{center}
\end{theorem}

Since there are $8$ primitive reduced forms,
$$
\fbox{$h(-840) = 8$,}
$$
or, simply, \textit{the class number is~$8$}.


\subsection{Kronecker's limit formula}

Let $(A,B,C)$ denote the positive definite quadratic form
$$
Q(X,Y) := A X^2 + 2B XY + C Y^2.
$$
We call
$$
m := AC - B^2
$$
the \textbf{determinant} of $Q(X,Y)$. These so-called \emph{Gauss
forms} are a special case of the general binary forms in the first
part, but we will only deal with these Gauss forms from now on. We
also point out that the \textit{determinant} of a Gauss form plays the
role of the \textit{discriminant} for the general binary forms. We
could express all of our results in terms of discriminants, instead of
determinants, but there is no need to do so.

The sum
$$
\fbox{$\dis S(s) := \sump\frac{1}{(AX^2 + 2BXY + CY^2)^s}\,,$}
$$
which is taken over all positive and negative integers $X$ and $Y$
with the exception of $X = Y = 0$, converges absolutely (and therefore
independently of the order of the terms) as long as $s > 1$, to a
finite value which is a \textit{function} of~$s$.

But, when $s \to 1^+$, the sum $S(s)$ \textit{diverges} to infinity,
and the Kronecker \textit{``Grenzformel''} gives us the first two
terms in the Laurent expansion of $S(s)$ around $s = 1$:
$$
S(s) = \frac{a_{-1}}{s-1} + a_0 + a_1(s - 1) +\cdots
$$

\begin{theorem}[Kronecker's Limit Formula]
$$
\fbox{$\dis \begin{aligned}
\lim_{s\to 1^+} &\biggl\{ \sump
\frac{1}{(AX^2 + 2BXY + CY^2)^s}
- \frac{\pi}{\sqrt{m}} \cdot \frac{1}{s-1} \biggr\}
\\[2\jot]
&\quad = -\frac{2\pi\gamma}{\sqrt{m}}
+ \frac{\pi}{\sqrt{m}} \ln\left(\frac{A}{4m}\right)
- \frac{2\pi}{\sqrt{m}} \ln \eta(\omega_1) \eta(\omega_2),
\end{aligned}$}
$$
\textit{where}
\begin{align}
\gamma &:= \text{ Euler's constant } = 0.5772\dots
\nonumber \\
\omega_1 &:= \frac{B + i\sqrt{m}}{A}, \qquad
\omega_2  := \frac{-B + i\sqrt{m}}{A},
\nonumber \\
\eta(\omega) &:= e^{\frac{\pi i\omega}{12}}
\prod_{n=1}^{\infty} (1 - e^{2\pi i n\omega}).
\tag*{\qed}
\end{align}
\end{theorem}

This statement of the Kronecker \textit{Grenzformel} is taken from
Weber~\cite{w1}. Kronecker uses his formula to prove the following
result.

\begin{theorem}[Fundamental Lemma]
If $(A,B,C)$ and $(A_1,B_1,C_1)$ are two forms with determinant $m$,
then
$$
\framebox{$ \begin{aligned}
\lim_{s\to 1^+} &\dis \biggl\{ \sump \frac{1}{(AX^2 + 2BXY + CY^2)^s}
- \sump \frac{1}{(A_1X^2 + 2B_1XY + C_1Y^2)^s} \biggr\}
\\[2\jot]
&= \dis \frac{2\pi}{\sqrt{m}}
\ln \biggl\{ \biggl(\frac{A}{A_1}\biggr)^{\frac{1}{2}} \cdot
\frac{\eta(\Omega_1)}{\eta(\omega_1)} \cdot
\frac{\eta(\Omega_2)}{\eta(\omega_2)} \biggr\}
\end{aligned} $}
$$
where
$$
\begin{aligned}
\omega_1 &:= \frac{B + i\sqrt{m}}{A},  &
\omega_2 &:= \frac{-B + i\sqrt{m}}{A},
\\
\Omega_1 &:= \frac{B_1 + i\sqrt{m}}{A_1},  &
\Omega_2 &:= \frac{-B_1 + i\sqrt{m}}{A_1}.
\end{aligned}
$$
\end{theorem}

\begin{proof}
Apply the \textit{Grenzformel} to $(A,B,C)$ and $(A_1,B_1,C_1)$
separately and subtract. The left hand side of the difference is the
left hand side of the stated proposition, while the right hand side
equals
\begin{align}
& \frac{\pi}{\sqrt{m}} \biggl\{ \ln\biggl(\frac{A}{4m}\biggr)
- \ln\biggl(\frac{A_1}{4m}\biggr)
- 2\bigl[\ln\eta(\omega_1)\eta(\omega_2)
- \ln\eta(\Omega_1)\eta(\Omega_2) \bigr] \biggr\}
\nonumber \\
&\qquad = \frac{\pi}{\sqrt{m}}
\biggl\{ \ln\biggl(\frac{A}{A_1}\biggr)^{\frac{1}{2}}
\biggl[\frac{\eta(\Omega_1)}{\eta(\omega_1)}\biggr]^2
\biggl[\frac{\eta(\Omega_2)}{\eta(\omega_2)}\biggr]^2 \biggr\}
\nonumber \\
&\qquad = \frac{2\pi}{\sqrt{m}}
\biggl\{ \ln\biggl(\frac{A}{A_1}\biggr)^{\frac{1}{2}}
\biggl[\frac{\eta(\Omega_1)}{\eta(\omega_1)}\biggr]
\biggl[\frac{\eta(\Omega_2)}{\eta(\omega_2)}\biggr] \biggr\} .
\tag*{\qed}
\end{align}
\hideqed
\end{proof}

\noindent
\textbf{Note:}\quad
$\omega_1$ and $\omega_2$ are called the \textbf{roots} of the form
$(A,B,C)$.


\subsection{The relation with Ramanujan's function \boldmath$g_n$}

The key observation which permits us to apply Kronecker's Fundamental
Lemma to the computation of $g_{210}$ is that the eight reduced forms
of determinant $210$ (or of discriminant $-840$) can be \textit{paired
off} as follows:
\begin{align*}
X^2 + 210Y^2 &\leftrightarrows 2X^2 + 105Y^2,
\\
3X^2 + 70Y^2 &\leftrightarrows 6X^2 + 35Y^2,
\\
5X^2 + 42Y^2 &\leftrightarrows 10X^2 + 21Y^2,
\\
7X^2 + 30Y^2 &\leftrightarrows 14X^2 + 15Y^2.
\end{align*}
The ``pairing'' is of the \textit{forms}
$$
\fbox{$AX^2 + 2CY^2 \ \leftrightarrows \ 2AX^2 + CY^2$}
$$
and the ``roots'' of the forms on each side are respectively $\omega$
and~$\Omega$, i.e.,
$$
\omega = \omega_1 = \omega_2 = \frac{\sqrt{-m}}{A}\ \leftrightarrows \
\Omega = \Omega_1 = \Omega_2 = \frac{\sqrt{-m}}{2A} = \frac{\omega}{2}
$$
that is, \textit{the root $\omega$ is paired off with the root
$\omega/2$}. Applying the Fundamental Lemma to this pair of forms
$(A,0,2C)$ and $(2A,0,C)$, we obtain
\begin{align*}
\lim_{s \to 1^+} &\dis \biggl\{ \sump \frac{1}{(AX^2 + 2CY^2)^s}
- \sump \frac{1}{(2AX^2 + CY^2)^s} \biggr\}
\\
&= \frac{2\pi}{\sqrt{m}} \ln \left\{ 
\biggl(\frac{A}{2A}\biggr)^{\frac{1}{2}} \,
\frac{\eta(\omega/2)}{\eta(\omega)} \,
\frac{\eta(\omega/2)}{\eta(\omega)} \right\}
\\
&= \frac{2\pi}{\sqrt{m}} \ln \biggl(\frac{1}{2}\biggr)^{\frac{1}{2}}
\biggl[\frac{\eta(\omega/2)}{\eta(\omega)}\biggr]^2 .
\end{align*}
But, if $q := e^{\pi i \omega}$, then
$$
\eta(\omega) = q^{\frac{1}{12}} \prod_{n=1}^\infty (1 - q^{2n}),
$$
and so
$$
\eta(\omega/2) = q^{\frac{1}{24}} \prod_{n=1}^\infty (1 - q^n).
$$
Therefore,
\begin{align*}
\frac{\eta(\omega/2)}{\eta(\omega)}
&= \frac{\dis q^{\frac{1}{24}} \prod_{n=1}^\infty (1 - q^n)}
{\dis q^{\frac{1}{12}} \prod_{n=1}^\infty (1 - q^{2n})}
\\
&= \frac{\dis q^{-\frac{1}{24}}
\prod_{n=1}^\infty (1 - q^{2n})(1 - q^{2n-1})}
{\dis \prod_{n=1}^\infty (1 - q^{2n})}
\\
&= q^{-\frac{1}{24}} \dis \prod_{n=1}^\infty (1 - q^{2n-1}).
\end{align*}
Since $B = 0$ for both of our forms, we conclude that
$\dis\omega = \frac{i\sqrt{m}}{A}$, and so the right-hand side becomes
$$
e^{-\frac{\pi\sqrt{m}}{24A}} \prod_{n=1}^\infty
\bigl[1 - e^{-(2n-1)\pi\sqrt{m}/A}\ \bigr]
= 2^{\frac{1}{4}}\, g_{m/A^2}
$$
as we see from the definition of $g_n$, given in equation
\eqref{eq:Rama-prod}. Therefore
\begin{align*}
\dis \lim_{s \to 1^+} \biggl\{ \sump \frac{1}{(AX^2 + 2CY^2)^s}
- \sump \frac{1}{(2AX^2 + CY^2)^s} \biggr\}
&= \frac{2\pi}{\sqrt{m}} \ln \biggl(\frac{1}{2}\biggr)^{\frac{1}{2}}
\bigl\{ 2^{\frac{1}{4}}\, g_{m/A^2} \bigr\}^2
\\
&= \frac{4\pi}{\sqrt{m}} \ln g_{m/A^2} \,.
\end{align*}

We have therefore obtained the following formula.

\begin{theorem}[Kronecker's Formula for {\boldmath $g_n$}]
$$
\lim_{s\to 1^+} \biggl\{ \sump \frac{1}{(AX^2 + 2CY^2)^s}
- \sump \frac{1}{(2AX^2 + CY^2)^s} \biggr\}
= \frac{4\pi}{\sqrt{m}}\, \ln g_{m/A^2} \,.
\eqno \mathbf{(G)}
$$
\end{theorem}

\textit{This is the formula which we will use to compute $g_{210}$}.

\vspace{6pt}

Kronecker, himself, expressed the right hand side of his
\textit{Grenzformel} in terms of the theta functions, while we have
chosen Ramanujan's notation.
This statement of the Kronecker formula in terms of the Ramanujan
function $g_n$ has not been explicitly stated in the literature
(although it is clearly implicitly contained in Weber's work).

Following earlier work of \textsc{Dirichlet}, Kronecker and
Weber~\cite{w1,w2} \textit{summed} both sides of \textbf{(G)} over all
reduced forms of determinant~$m$. However, and here comes another
``brilliancy,'' they \textit{weighted} the sum by multiplying each of
the $h(m) \equiv h(\Delta)$ summands by the Jacobi symbol
$\dis\biggl( \frac{\delta}{A + 2C} \biggr)$ (see Niven and Zuckerman
\cite{nz}) where $\delta$ is a \textit{fundamental discriminant}. We
will define this term later on.

Since the weights are numerically equal to $+1$ or $-1$, one hopes for
useful \textit{cancellations}. Miraculously \emph{everything} cancels
(!) and one obtains a sum computable from the Dirichlet \textit{class
number formulas}.


\subsection{The Kronecker--Weber computation of the weighted sum for
\boldmath $m = 210$}

The sum which we wish to compute is

\noindent
\begin{align*}
R := \lim_{s \to 1^+} \Biggl\{ 
&\dis \biggl(\frac{\delta}{1}\biggr)
\biggl[ \sump \frac{1}{(X^2 + 210Y^2)^s}
- \sump \frac{1}{(2X^2 + 105Y^2)^s} \biggr]
\\[\jot]
&\dis + \biggl(\frac{\delta}{107}\biggr)
\biggl[ \sump \frac{1}{(2X^2 + 105Y^2)^s}
- \sump \frac{1}{(X^2 + 210Y^2)^s} \biggr]
\\[\jot]
&\dis + \biggl(\frac{\delta}{73}\biggr)
\biggl[ \sump \frac{1}{(3X^2 + 70Y^2)^s}
- \sump \frac{1}{(6X^2 + 35Y^2)^s} \biggr]
\\[\jot]
&\dis + \biggl(\frac{\delta}{47}\biggr)
\biggl[ \sump \frac{1}{(5X^2 + 42Y^2)^s}
- \sump \frac{1}{(10X^2 + 21Y^2)^s} \biggr]
\\[\jot]
&\dis + \biggl(\frac{\delta}{41}\biggr)
\biggl[ \sump \frac{1}{(6X^2 + 35Y^2)^s}
- \sump \frac{1}{(3X^2 + 70Y^2)^s} \biggr]
\\[\jot]
&\dis + \biggl(\frac{\delta}{37}\biggr)
\biggl[ \sump \frac{1}{(7X^2 + 30Y^2)^s}
- \sump \frac{1}{(14X^2 + 15Y^2)^s} \biggr]
\\[\jot]
&\dis + \biggl(\frac{\delta}{31}\biggr)
\biggl[ \sump \frac{1}{(10X^2 + 21Y^2)^s}
- \sump \frac{1}{(5X^2 + 42Y^2)^s} \biggr]
\\[\jot]
&\dis + \biggl(\frac{\delta}{29}\biggr)
\biggl[ \sump \frac{1}{(14X^2 + 15Y^2)^s}
- \sump \frac{1}{(7X^2 + 30Y^2)^s} \biggr] \Biggr\} \,,
\end{align*}
where $\delta$ is a \textbf{fundamental discriminant}. By definition,
this means that
\begin{enumerate}
\item $\delta$ divides $210$, and
\item $\delta \equiv 1 \pmod 4$,
\end{enumerate}
which implies that
$$
\delta = 1, \ -3, \ 5, \ -7,\ -15, \ 21, \ -35, \ 105.
$$

Kronecker and Weber compute this sum in two different ways. Following
their lead, \textit{first} we will use Kronecker's formula \textbf{(G)}
and obtain linear combinations of Ramanujan's function~$g_n$.
\textit{Second}, we will use the Dirichlet class number formulas.

To carry out the first computation, we compute the following
table of Jacobi symbols.
\begin{center}
\begin{tabular}{|c||c|c|c|c|c|c|c|c|}
\hline
$\delta=$ & 1 & $-3$ & 5 & $-7$ & $-15$ & 21 & $-35$ & 105
\\ \hline\hline
\rule[-3mm]{0mm}{8mm} $(\frac{\delta}{1})$
& 1 & 1 & 1 & 1 & 1 & 1 & 1 & 1
\\ \hline
\rule[-3mm]{0mm}{8mm}$(\frac{\delta}{107})$
& 1 & $-1$ & $-1$ & 1 & 1 & $-1$ & $-1$ & 1
\\ \hline
\rule[-3mm]{0mm}{8mm}$(\frac{\delta}{73})$
& 1 & 1 & $-1$ & $-1$ & $-1$ & $-1$ & 1 & 1
\\ \hline
\rule[-3mm]{0mm}{8mm}$(\frac{\delta}{47})$
& 1 & $-1$ & $-1$ & $-1$ & 1 & 1 & 1 & $-1$
\\ \hline
\rule[-3mm]{0mm}{8mm}$(\frac{\delta}{41})$
& 1 & $-1$ & 1 & $-1$ & $-1$ & 1 & $-1$ & 1
\\ \hline
\rule[-3mm]{0mm}{8mm}$(\frac{\delta}{37})$
& 1 & 1 & $-1$ & 1 & $-1$ & 1 & $-1$ & $-1$
\\ \hline
\rule[-3mm]{0mm}{8mm}$(\frac{\delta}{31})$
& 1 & 1 & 1 & $-1$ & 1 & $-1$ & $-1$ & $-1$
\\ \hline
  \rule[-3mm]{0mm}{8mm}$(\frac{\delta}{29})$
& 1 & $-1$ & 1 & 1 & $-1$ & $-1$ & 1 & $-1$
\\ \hline
\end{tabular}
\rule[-3mm]{0mm}{8mm}
\end{center}

Applying the Kronecker formula (G) to the ``pairs'', we find:
\begin{align*}
\dis \lim_{s \to 1^+} \biggl[ 
\sump \frac{1}{(X^2 + 210Y^2)^s}
- \sump \frac{1}{(2X^2 + 105Y^2)^s} \biggr]
&= \dfrac{4\pi}{\sqrt{210}}\, \ln g_{210} \,,
\\
\dis \lim_{s \to 1^+} \biggl[
\sump \frac{1}{(3X^2 + 70Y^2)^s}
- \sump \frac{1}{(6X^2 + 35Y^2)^s} \biggr]
&= \dfrac{4\pi}{\sqrt{210}}\, \ln g_{210/3^2} \,,
\\
\dis \lim_{s \to 1^+} \biggl[
\sump \frac{1}{(5X^2 + 42Y^2)^s}
- \sump \frac{1}{(10X^2 + 21Y^2)^s} \biggr]
&= \dfrac{4\pi}{\sqrt{210}}\, \ln g_{210/5^2} \,,
\\
\dis \lim_{s \to 1^+}\biggl[
\sump \frac{1}{(14X^2 + 15Y^2)^s}
- \sump \frac{1}{(7X^2 + 30Y^2)^s} \biggr]
&= \dfrac{4\pi}{\sqrt{210}}\, \ln g_{210/7^2} \,.
\end{align*}
Therefore, the sum $R$ becomes
\begin{align}
R &= \biggl[ \biggl(\dfrac{\delta}{1}\biggr)
- \biggl(\dfrac{\delta}{107}\biggr) \biggr]
\dfrac{4\pi}{\sqrt{210}}\ln g_{210}
+ \biggl[ \biggl(\dfrac{\delta}{73}\biggr)
- \biggl(\dfrac{\delta}{41}\biggr) \biggr]
\dfrac{4\pi}{\sqrt{210}}\ln g_{210/3^2}
\nonumber\\
&\qquad + \biggl[ \biggl(\dfrac{\delta}{47}\biggr)
- \biggl(\dfrac{\delta}{31}\biggr) \biggr]
\dfrac{4\pi}{\sqrt{210}}\ln g_{210/5^2}
+ \biggl[ \biggl(\dfrac{\delta}{37}\biggr)
- \biggl(\dfrac{\delta}{29}\biggr) \biggr]
\dfrac{4\pi}{\sqrt{210}}\ln g_{210/7^2} \,.
\label{eq:R-sum}
\end{align}
{}From the table, we see that
\begin{center}
\begin{tabular}{|c|c|c|c|c|}
\hline
\rule[-3mm]{0mm}{8mm}$\delta$ &
$\biggl(\dfrac{\delta}{1}\biggr) - \biggl(\dfrac{\delta}{107}\biggr)$ &
$\biggl(\dfrac{\delta}{73}\biggr) - \biggl(\dfrac{\delta}{41}\biggr)$ &
$\biggl(\dfrac{\delta}{47}\biggr) - \biggl(\dfrac{\delta}{31}\biggr)$ &
$\biggl(\dfrac{\delta}{37}\biggr) - \biggl(\dfrac{\delta}{29}\biggr)$
\\ \hline\hline
1 & 0 & 0 & 0 & 0 \\ \hline
$-3$ & 2 & 2 & $-2$ & 2 \\ \hline
5 & 2 & $-2$ & $-2$ & $-2$ \\ \hline
$-7$ & 0 & 0 & 0 & 0 \\ \hline
$-15$ & 0& 0& 0& 0 \\ \hline
21 & 2 & $-2$ & 2 & 2 \\ \hline
$-35$ & 2 & 2 & 2 & $-2$ \\ \hline
105 & 0 & 0 & 0 & 0 \\ \hline
\hline
\end{tabular}
\end{center}
and four sums survive:
\begin{center}
\fbox{\textbf{Surviving Sums}}
\end{center}
\begin{center}
\begin{tabular}{|c||c|}
\hline
$\delta$ & \textbf{SUM} \\
\hline\hline
\rule[-4mm]{0mm}{12mm}$-3$ & $\dfrac{8\pi}{\sqrt{210}}
\left\{\ln g_{210} + \ln g_{210/3^2} - \ln g_{210/5^2}
+ \ln g_{210/7^2} \right\}$
\\ \hline
\rule[-4mm]{0mm}{12mm}5 & $\dfrac{8\pi}{\sqrt{210}}
\left\{\ln g_{210} - \ln g_{210/3^2} - \ln g_{210/5^2}
- \ln g_{210/7^2} \right\}$
\\ \hline
\rule[-4mm]{0mm}{12mm}21 & $\dfrac{8\pi}{\sqrt{210}}
\left\{\ln g_{210} - \ln g_{210/3^2} + \ln g_{210/5^2}
+ \ln g_{210/7^2} \right\}$
\\ \hline
\rule[-4mm]{0mm}{12mm}$-35$ & $\dfrac{8\pi}{\sqrt{210}}
\left\{\ln g_{210} + \ln g_{210/3^2} + \ln g_{210/5^2}
- \ln g_{210/7^2} \right\}$
\\ \hline
\end{tabular}
\end{center}

If we sum these four survivors, the \textit{miraculous cancellation
takes place} and we obtain:
\begin{equation}
\fbox{$\dis R = \frac{8\pi}{\sqrt{210}} \cdot 4 \cdot \ln g_{210}
= \frac{32\pi}{\sqrt{210}}\, \ln g_{210} $}
\label{eq:miracle-R}
\end{equation}
a truly astonishing result!

This method of summing the series R to obtain the four series in
\eqref{eq:R-sum} to obtain the final sum \eqref{eq:miracle-R} taken
over all values of the remaining $\delta$'s is Kronecker's novel
contribution. Weber \cite{w1,w2} refined it into a broadly applicable
tool. We have presented the simplest version possible so as to
recapture the classical beauty and elegance of Kronecker's original
version.


\subsection{The Dirichlet computation of the weighted sum}

We group the coefficients of the \textit{same} summand in $R$ so
that $R$ becomes:

\noindent
\begin{align*}
R = \lim_{s\to 1^+} \Biggl\{ 
&\dis \biggl[ \biggl(\frac{\delta}{1}\biggr)
- \biggl(\frac{\delta}{107}\biggr) \biggr]
\sump \frac{1}{(X^2 + 210Y^2)^s}
\\[\jot]
&\dis + \biggl[ \biggl(\frac{\delta}{107}\biggr)
- \biggl(\frac{\delta}{1}\biggr) \biggr]
\sump \frac{1}{(2X^2 + 105Y^2)^s}
\\[\jot]
&\dis + \biggl[ \biggl(\frac{\delta}{73}\biggr)
- \biggl(\frac{\delta}{41}\biggr) \biggr]
\sump \frac{1}{(3X^2 + 70Y^2)^s}
\\[\jot]
&\dis + \biggl[ \biggl(\frac{\delta}{47}\biggr)
- \biggl(\frac{\delta}{31}\biggr) \biggr]
\sump \frac{1}{(5X^2 + 42Y^2)^s}
\\[\jot]
&\dis + \biggl[ \biggl(\frac{\delta}{41}\biggr)
- \biggl(\frac{\delta}{73}\biggr) \biggr]
\sump \frac{1}{(6X^2 + 35Y^2)^s}
\\[\jot]
&\dis + \biggl[ \biggl(\frac{\delta}{37}\biggr)
- \biggl(\frac{\delta}{29}\biggr) \biggr]
\sump \frac{1}{(7X^2 + 42Y^2)^s}
\\[\jot]
&\dis + \biggl[ \biggl(\frac{\delta}{31}\biggr)
- \biggl(\frac{\delta}{47}\biggr) \biggr]
\sump \frac{1}{(10X^2 + 21Y^2)^s}
\\[\jot]
&\dis + \biggl[ \biggl(\frac{\delta}{29}\biggr)
- \biggl(\frac{\delta}{37}\biggr) \biggr]
\sump \frac{1}{(14X^2 + 15Y^2)^s} \Biggr\} \,.
\end{align*}

To see what's going on, we examine a particular sum, say
\begin{equation}
\biggl[ \biggl(\frac{\delta}{41}\biggr)
- \biggl(\frac{\delta}{73}\biggr) \biggr]
\sump \frac{1}{(6X^2 + 35Y^2)^s} \,.
\label{eq:fifth-sum}
\end{equation}
First we observe that
$$
\biggl(\frac{\delta}{41}\biggr) = - \biggl(\frac{\delta}{73}\biggr),
$$
so that \eqref{eq:fifth-sum} can be written
$$
2 \biggl(\frac{\delta}{41}\biggr) \sump \frac{1}{(6X^2 + 35Y^2)^s} \,.
$$
Let's look at a particular summand.  For example, take
$$
X = 3 \, Y = 7.
$$
Then $6 X^2 + 35 Y^2 = 1769$. The summand $\dfrac{1}{(1769)^s}$
appears once for every pair of integers $(X,Y)$ which satisfies
$$
6 X^2 + 35 Y^2 = 1769.
$$
Since, as it is
easy to verify (see below), the set of solutions is
$$
(X,Y) = (\pm 3,\pm 7)  \phrase{and}  (\pm 17,\pm 1),
$$
each of the possible eight pairs contributes $1$ to the coefficient of
$\dfrac{1}{(1769)^s}$. We therefore conjecture that \textit{the total
contribution of the denominator $(1769)^s$ to the sum $R$ is equal to}
$$
\frac{8}{(1769)^s} \,.
$$
How do we \textit{know} that the eight listed pairs $(X,Y)$ exhaust
all solutions of $6X^2 + 35Y^2 = 1769$, or that \textit{there are no
other summands from the other reduced forms of determinant $210$}?

Dirichlet, himself, answered this question when he found \textit{the
formula for the total number of proper representations of an
integer~$N$ by the forms of determinant~$m$}: see Cox~\cite{cox}. If
$f(N)$ denotes this number of representations, then
$$
\fbox{$\dis 
f(N) = 2 \sum_{d|N,\; d\geq 1} \biggl(\frac{-m}{d}\biggr) $}
$$
where $\dis \biggl(\frac{-m}{d}\biggr)$ is the Jacobi symbol. In
our case
$$
f(1729) = 2\left\{\biggl( \frac{-210}{1}\biggr) + \biggl(
\frac{-210}{29}\biggr) + \biggl( \frac{-210}{61}\biggr) + \biggl(
\frac{-210}{61}\biggr)\right\} = 8,
$$
which verifies that the eight pairs we displayed are \textit{all} of
them.

Thus, the total contribution to the sum $R$ is
\begin{equation}
\frac{\dis 4 \biggl(\frac{\delta}{41}\biggr)\left\{\biggl(
\frac{-210}{1}\biggr) + \biggl( \frac{-210}{29}\biggr) + \biggl(
\frac{-210}{61}\biggr) + \biggl(
\frac{-210}{1769}\biggr)\right\}}{(1769)^s} \,.
\label{eq:R-total}
\end{equation}
Now we use some properties of the Jacobi symbol:
$$
\biggl(\frac{m}{p}\biggr) \, \biggl(\frac{m}{p}\biggr)
= \dis \biggl(\frac{mn}{p}\biggr)  \phrase{and} 
\biggl(\frac{m}{pq}\biggr) \, \biggl(\frac{m}{q}\biggr)
= \dis \biggl(\frac{m}{p}\biggr),
$$
as well as
$$
\biggl(\frac{\delta}{41}\biggr) = \biggl(\frac{\delta}{1769}\biggr).
$$
Therefore, if $\delta\delta ' = 210$, the expression in
\eqref{eq:R-total} equals
\begin{align*}
& \frac{\dis 4\biggl(\frac{\delta}{41}\biggr)\left\{\biggl(
\frac{-210}{1}\biggr) + \biggl( \frac{-210}{29}\biggr) + \biggl(
\frac{-210}{61}\biggr) + \biggl(
\frac{-210}{1769}\biggr)\right\}}{(1769)^s}
\\[\jot]
&= \frac{\dis 4\biggl(\frac{\delta}{1769}\biggr)\left\{\biggl(
\frac{-210}{1}\biggr) + \biggl( \frac{-210}{29}\biggr) + \biggl(
\frac{-210}{61}\biggr) + \biggl(
\frac{-210}{1769}\biggr)\right\}}{(1769)^s}
\\[\jot]
&= \frac{\dis 4\biggl(\frac{\delta}{1769}\biggr)
\left\{\biggl( \frac{\delta }{1}\biggr)\biggl( \frac{\delta '}{1}\biggr)
 + \biggl( \frac{\delta }{29}\biggr)\biggl( \frac{\delta '}{29}\biggr)
 + \biggl( \frac{\delta }{61}\biggr)\biggl( \frac{\delta '}{61}\biggr) + 
\biggl( \frac{\delta }{1769}\biggr)\biggl( \frac{\delta
'}{1769}\biggr)\right\}}{(1769)^s}
\\[\jot]
&= \frac{\dis 4\left\{\biggl( \frac{\delta}{1769}\biggr)
 \biggl( \frac{\delta'}{1}\biggr)
 + \biggl( \frac{\delta }{61}\biggr)\biggl( \frac{\delta '}{29}\biggr)
 + \biggl( \frac{\delta }{29}\biggr)\biggl( \frac{\delta '}{61}\biggr)
 + \biggl( \frac{\delta }{1}\biggr)\biggl( \frac{\delta'}{1769}\biggr)
\right\}}{(1769)^s} \,.
\end{align*}
But this is the coefficient of $\dfrac{1}{(1769)^s}$ in the expansion
of the \textit{product series}:
\begin{equation}
4 \Biggl\{ \sum_{n=1}^\infty
\biggl(\frac{\delta}{n}\biggr) \frac{1}{n^s} \Biggr\}
\Biggl\{ \sum_{m=1}^\infty 
\biggl( \frac{\delta}{m}\biggr) \frac{1}{m^s} \Biggr\} \,,
\label{eq:prod-series}
\end{equation}
and we conclude that
\begin{center}
\fbox{\textit{The sum $R$ is equal to the product
\eqref{eq:prod-series}}.}
\end{center}
We boxed this statement since it is, in fact, a \textbf{fundamental
theorem} in the theory of these series.

Therefore,  the first line in our table of surviving sums is
\begin{align*}
\frac{8\pi}{\sqrt{210}} 
& \left\{\ln g_{210} + \ln g_{210/3^2}
- \ln g_{210/5^2} + \ln g_{210/7^2} \right\}
\\
&= \dis 4\left\{\sum_{n=1}^{\infty}
\biggl( \frac{-3}{n}\biggr)\frac{1}{n^s} \right\}
\left\{\sum_{m=1}^{\infty}\biggl( \frac{70 }{m}\biggr)\frac{1}{m^s}
\right\}
\\
&= \dis 4\left\{\sum_{n=1}^{\infty}
\biggl( \frac{-3}{n}\biggr)\frac{1}{n^s} \right\}
\left\{ \sum_{m=1}^{\infty}
\biggl(\frac{280 }{m}\biggr)\frac{1}{m^s} \right\} \,.
\end{align*}

If we carry out the same computation for the other three lines of
surviving sums we obtain the following four equations:
\begin{align}
\dfrac{8\pi}{\sqrt{210}} \left\{
\ln g_{210} + \ln g_{210/3^2} - \ln g_{210/5^2} + \ln g_{210/7^2}
\right\}
&= \dis 4\left\{\sum_{n=1}^{\infty}\biggl(
\frac{-3}{n}\biggr)\frac{1}{n^s} \right\}
\left\{\sum_{m=1}^{\infty}\biggl( \frac{280 }{m}\biggr)\frac{1}{m^s}
\right\} \,,
\nonumber \\[2\jot]
\dfrac{8\pi}{\sqrt{210}} \left\{
\ln g_{210} - \ln g_{210/3^2} - \ln g_{210/5^2} - \ln g_{210/7^2}
\right\}
&= \dis 4\left\{\sum_{n=1}^{\infty}\biggl(
\frac{5}{n}\biggr)\frac{1}{n^s} \right\}
\left\{\sum_{m=1}^{\infty}\biggl( \frac{-168}{m}\biggr)\frac{1}{m^s}
\right\} \,,
\nonumber \\[2\jot]
\dfrac{8\pi}{\sqrt{210}} \left\{
\ln g_{210} - \ln g_{210/3^2} + \ln g_{210/5^2} + \ln g_{210/7^2}
\right\}
&= \dis 4\left\{\sum_{n=1}^{\infty}\biggl(
\frac{21}{n}\biggr)\frac{1}{n^s} \right\}
\left\{\sum_{m=1}^{\infty}\biggl( \frac{-40}{m}\biggr)\frac{1}{m^s}
\right\} \,,
\nonumber \\[2\jot]
\dfrac{8\pi}{\sqrt{210}} \left\{
\ln g_{210} + \ln g_{210/3^2} + \ln g_{210/5^2} - \ln g_{210/7^2}
\right\}
&= \dis 4 \left\{ \sum_{n=1}^{\infty}\biggl(
\frac{-35}{n}\biggr)\frac{1}{n^s} \right\}
\left\{\sum_{m=1}^{\infty}\biggl( \frac{24}{m}\biggr)\frac{1}{m^s}
\right\} \,.
\label{eq:four-eqns}
\end{align}

\subsubsection{The Dirichlet class number formulas}

Onw of the glories of 19th century mathematics is Dirichlet's
formulas for the limit of the sums of the series appearing in
\eqref{eq:four-eqns} as $s \to 1^+$. Namely,
\begin{align}
\dis \sum_{n=1}^{\infty}\biggl(
\frac{\delta}{n}\biggr)\frac{1}{n} 
&= \dis\frac{\pi}{\sqrt{-\delta}} \cdot K(\delta)
\phrase{if} \delta < 0,
\label{eq:Diri-formula-neg}
\\[2\jot]
\dis \sum_{n=1}^{\infty}\biggl(
\frac{\delta}{n}\biggr)\frac{1}{n} 
&= \dis\frac{\ln\varepsilon}{\sqrt{\delta}} \cdot K(\delta)
\phrase{if} \delta > 0,
\label{eq:Diri-formula-pos}
\end{align}
where $K(\delta)$ is the number of properly
primitive classes of discriminant $\delta$ and
$$
\varepsilon := \dfrac{T + U\sqrt{\delta}}{2}
$$
is the \textit{minimal solution of the ``even'' Pell equation}:
$$
T^2 - \delta U^2 = 4.
$$
These formulas are developed in numerous books including
Weber~\cite{w1}.

\subsubsection{The final computation of \boldmath $g_{210}$}

We make a table:
\begin{center}
\begin{tabular}{|c|c|c|c|c|}
\hline
\rule[-4mm]{0mm}{12mm} $\delta$ & $K(\delta)$ &$ \delta '$ &$
K(\delta ')$ & $\dis \frac{T + U\sqrt{\delta}}{2}$
\\\hline\hline
\rule[-4mm]{0mm}{12mm} $-3$ & $\dis \frac{1}{3}$ &$280$ &$4$ &
$\dis \frac{502 + 30\sqrt{280}}{2}$
\\\hline
\rule[-4mm]{0mm}{12mm} $5$ & $1$ & $-168$& $4$ & $\dis
\frac{3 + \sqrt{5}}{2}$
\\\hline
\rule[-4mm]{0mm}{12mm} $21$ & $2$ & $-40$ & $2$&  $\dis
\frac{5 + \sqrt{21}}{2}$
\\\hline
\rule[-4mm]{0mm}{12mm} $-35$ & $2$ & $24$ & 2&
$\dis\frac{10 + 2\sqrt{24}}{2}$\\
\hline
\end{tabular}
\end{center}
\begin{center}
\fbox{\textbf{Class Numbers and Units}}
\end{center}

The values $K(\delta)$ in the table can be found in
Flath~\cite{flath}. Substituting these numbers in
\eqref{eq:Diri-formula-neg} and \eqref{eq:Diri-formula-pos} we obtain,
for example,
$$
\sum_{n=1}^{\infty} \biggl( \frac{-3}{n}\biggr) \frac{1}{n}
= \frac{\pi}{\sqrt{3}} \cdot K(-3)
= \frac{\pi}{3\sqrt{3}} \,,
$$
and
\begin{align*}
\sum_{n=1}^{\infty} \biggl( \frac{280 }{n}\biggr) \frac{1}{n}
&= \dis \frac{K(280)}{\sqrt{280}} \cdot
\ln \biggl( \frac{T + U\sqrt{280}}{2} \biggr)
\\
&= \dis \frac{4}{\sqrt{280}}\  \cdot
 \ln \biggl(\frac{502 + 30\sqrt{280}}{2} \biggr)
\\
&= \dis \frac{4}{\sqrt{280}}\  \cdot
 \ln \bigl(5\sqrt{5} + 3\sqrt{14}\bigr)^2
\\
&= \dis \frac{8}{\sqrt{280}}\  \cdot
 \ln \bigl(5\sqrt{5} + 3\sqrt{14}\bigr).
\end{align*}
Both of these summations have independent interest.

Therefore,
$$
\frac{2\pi}{\sqrt{210}}\bigl\{\ln g_{210} + \ln
g_{210/3^2} - \ln g_{210/5^2} + \ln
g_{210/7^2} \bigr\} = \frac{\pi}{3\sqrt{3}} \cdot 
\frac{8}{\sqrt{280}} \cdot \ln \bigl(5\sqrt{5} + 3\sqrt{14}\bigr),
$$
or
\begin{equation}
\ln g_{210} + \ln g_{210/3^2} - \ln g_{210/5^2} + \ln g_{210/7^2}
=  \frac{1}{3} \ln \bigl(5\sqrt{5} + 3\sqrt{14}\bigr)^2 .
\label{eq:log-sum}
\end{equation}

The Dirichlet class number formulas applied the three remaining
product series give us
\begin{align*}
\dis 4\left\{\sum_{n=1}^{\infty}\biggl( \frac{5
}{n}\biggr)\frac{1}{n^s} \right\}\dis
\left\{\sum_{m=1}^{\infty}\biggl( \frac{-168}{m}\biggr)\frac{1}{m^s}
\right\}
&= \frac{\pi}{\sqrt{210}}\frac{1 \cdot 4}{2}\ln\biggl(
\frac{3 + \sqrt{5}}{2}\biggr) \cdot 4,
\\[\jot]
\dis 4\left\{\sum_{n=1}^{\infty}\biggl( \frac{21
}{n}\biggr)\frac{1}{n^s} \right\}\dis
\left\{\sum_{m=1}^{\infty}\biggl( \frac{-40}{m}\biggr)\frac{1}{m^s}
\right\}
&= \frac{\pi}{\sqrt{210}}\frac{2 \cdot 2}{2}\ln\biggl(
\frac{5 + \sqrt{21}}{2}\biggr) \cdot 4,
\\[\jot]
\dis 4\left\{\sum_{n=1}^{\infty}\biggl( \frac{-35
}{n}\biggr)\frac{1}{n^s} \right\}\dis
\left\{\sum_{m=1}^{\infty}\biggl( \frac{24}{m}\biggr)\frac{1}{m^s}
\right\}
&= \frac{\pi}{\sqrt{210}}\frac{2 \cdot 2}{2}\ln\biggl(
\frac{10 + 2\sqrt{24}}{2}\biggr) \cdot 4.
\end{align*}
Thus the last three lines of \eqref{eq:four-eqns} permit us to
conclude
\begin{align}
\ln g_{210} - \ln g_{210/3^2} - \ln g_{210/5^2} - \ln g_{210/7^2} 
&= \dis \ln\biggl( \frac{3 + \sqrt{5}}{2}\biggr),
\nonumber \\[2\jot]
\ln g_{210} - \ln g_{210/3^2} + \ln g_{210/5^2} + \ln g_{210/7^2} 
&= \dis \ln\biggl(\frac{5 + \sqrt{21}}{2}\biggr),
\nonumber \\[2\jot]
\ln g_{210} + \ln g_{210/3^2} + \ln g_{210/5^2} - \ln g_{210/7^2} 
&= \dis \ln \biggl(\frac{10 + 2\sqrt{24}}{2}\biggr).
\label{eq:log-sum-more}
\end{align}
Summing \eqref{eq:log-sum} and \eqref{eq:log-sum-more}, we obtain
\begin{align*}
4\,\ln g_{210}
&= \frac{2}{3}\ln\biggl(5\sqrt{5} + 3\sqrt{14}\biggr) + \ln\biggl(\dis
\frac{3 + \sqrt{5}}{2}\biggr) + \ln\biggl(\dis
\frac{5 + \sqrt{21}}{2}\biggr) + \ln \biggl(\dis
\frac{10 + 2\sqrt{24}}{2}\biggr)
\\
 &= \ln \left[
\biggl(5\sqrt{5} + 3\sqrt{14}\biggr)^{\frac{2}{3}} \cdot \biggl(\dis
\frac{3 + \sqrt{5}}{2}\biggr) \cdot \biggl(\dis
\frac{5 + \sqrt{21}}{2}\biggr) \cdot \biggl(\dis
\frac{10 + 2\sqrt{24}}{2}\biggr)\right]
\\
\implies \ g_{210} &= 
\biggl(5\sqrt{5} + 3\sqrt{14}\biggr)^{\frac{1}{6}} \cdot \biggl(\dis
\frac{3 + \sqrt{5}}{2}\biggr)^{\frac{1}{4}} \cdot \biggl(\dis
\frac{5 + \sqrt{21}}{2}\biggr)^{\frac{1}{4}} \cdot \biggl(\dis
\frac{10 + 2\sqrt{24}}{2}\biggr)^{\frac{1}{4}}
\\
 &= \biggl(5\sqrt{5} + 3\sqrt{14}\biggr)^{\frac{1}{6}} \cdot  \biggl(\dis
\frac{1 + \sqrt{5}}{2}\biggr)^{\frac{1}{2}} \cdot  \biggl(\dis
\frac{\sqrt{3} + \sqrt{7}}{2}\biggr)^{\frac{1}{2}} \cdot 
\biggl(\sqrt{2} + \sqrt{3}\biggr)^{\frac{1}{2}}
\end{align*}
That is,
{\large \boldmath $$
\fbox{$\dis
g_{210} = \sqrt{\sqrt{2} + \sqrt{3}} \cdot 
\sqrt[6]{\mathbf{5\sqrt{5} + 3\sqrt{14}}} \cdot
\sqrt{\frac{\sqrt{3} + \sqrt{7}}{2}} \cdot
\sqrt{\frac{\sqrt{5} + 1}{2}} $}
$$}
\emph{This \textbf{is} the value of $g_{210}$ which Watson
copied from Ramanujan's first draft of his notebook!}


\subsection{Final comments}

Weber showed that these beautiful and miraculous computations are
instances of general theorems, and \textit{always occur} for each of
the $65$ idoneal numbers, and not just for $m = 210$. His classic
exposition, however, never deals with the case of \textit{even} $m$
which (as we saw) converts the Kronecker Limit Formula into a
statement about the elliptic modular function
$\dfrac{\eta(\omega/2)}{\eta(\omega)}$. Instead, he develops the
entire theory for the related function
$\dis e^{-\pi i/24} \frac{\eta((\omega + 1)/2)}{\eta(\omega)}$. Since
the discussion of the former function is nowhere to be found, we will
carry it out in detail.


\subsection{The miracle explained: Kronecker's theory}

Weber \cite{w1,w2} showed that the ``miraculous cancellation'' is no
accident but rather is explainable \textit{a priori}, and we will now
develop his theory. Our proofs are modelled on Weber's, but since we
consider a very special case, we achieve considerable simplifications.

We make \textit{two} simplifying assumptions:
\begin{enumerate}
\item
\textit{All reduced forms of determinant $m$ are of the form}
$$
A X^2 + C Y^2.
$$
This means that the determinant $m = AC$ is a  \textit{convenient
number} \cite[Thm.~3.22(ii)]{cox}.
\item
$m = 2P$ where $P$ is a product of $t$ \textit{distinct odd primes}.
\end{enumerate}
Therefore (as is easy to see directly, or also is proven in
Theorem 3.22(v) of Cox~\cite{cox}),
$$
h(m) = 2^t,
$$
and the reduced forms break up into $2^{t-1}$ forms
$$
Q = A X^2 + 2C Y^2
$$
and $2^{t-1}$ homologues
$$
Q = 2A X^2 + C Y^2.
$$

Accordingly, we can \textit{pair them off}:
$$
\begin{aligned}
Q_1 := A_1 X^2 + 2C_1 Y^2
\quad &\leftrightarrows\quad
Q'_1 := 2A_1X^2 + C_1Y^2,
\\
Q_2 :=  A_2X^2 + 2C_2Y^2
\quad &\leftrightarrows\quad
Q'_2 := 2A_2X^2 + C_2Y^2,
\\
\vdots \qquad\qquad\qquad & \qquad\qquad\qquad \vdots
\\
Q_T :=  A_T X^2 + 2C_T Y^2
\quad &\leftrightarrows\quad
Q'_T := 2A_T X^2 + C_T Y^2,
\end{aligned}
$$
where
$$
T := 2^{t-1}.
$$

\begin{definition}
The number $\delta$ is a \textbf{fundamental discriminant} if and
only if
\begin{enumerate}
\item
$\delta \equiv 1 \pmod 4$ and $\delta$ is square-free; or
\item
$\delta  = 4\delta_1$, where $\delta_1 \not\equiv 1 \pmod 4$ and
$\delta_1$ is square-free.
\end{enumerate}
\end{definition}

Thus, \textit{every square-free odd number}, taken with \textit{appropriate sign},
is a fundamental discriminant.

The number of fundamental discriminants which divide our discriminant
$D = -4m$ is therefore \emph{the number of odd divisors of} $m$. But the
number of odd divisors of $m$ is equal to
$$
2^{t} = h(m).
$$

We associate the Jacobi symbol
$$
\chi(\delta;A,C) := \biggl(\frac{\delta}{A + C}\biggr)
\equiv \chi(\delta;Q)
$$
to the form $Q := A X^2 + B Y^2$.

Then the sum, $R$, which we wish to compute is
\begin{align*}
R(\delta) := R
&:= \dis\lim_{s\to 1^+} \Bigg\{ \chi(\delta;Q_1)
\left[\sump\frac{1}{Q_1^s} - \sump\frac{1}{Q_1'^s}\right]
+ \chi(\delta;Q_1')
\left[\sump \frac{1}{Q_1'^s} - \sump \frac{1}{Q_1^s}\right]
\\[\jot]
&\dis\qquad + \chi(\delta;Q_2)
\left[ \sump \frac{1}{Q_2^s} - \sump
 \frac{1}{Q_2'^s}\right]  + \chi(\delta;Q_2')
\left[ \sump \frac{1}{Q_2'^s} - \sump \frac{1}{Q_2^s}\right]
+ \cdots
\\[\jot]
&\dis\qquad + \chi(\delta;Q_{h})
\left[ \sump \frac{1}{Q_{h}^s} - \sump
 \frac{1}{Q_{h}'^s}\right]  + \chi(\delta;Q_{h}')
\left[ \sump \frac{1}{Q_{h}'^s} - \sump \frac{1}{Q_{h}^s}\right]
\Bigg\} \,,
\end{align*}
where $\delta$ is any one of the $2^{t}$ fundamental discriminants
which divide $D = - 4m$.

Our final result comes from computing
$$
H := \sum_{\delta} R(\delta),
$$
where $\delta$ runs over the fundamental discriminants which divide
$-4m$.

\subsubsection{The first cancellation: Kronecker's computation}

The key step in Kronecker's computation is the following Lemma which
shows what happens when we compute the contribution, $R_1$, of the two
``paired off'' forms, $Q$ and $Q'$, to the sum $R$.

\begin{lemma}
\label{lm:R-one}
\begin{align*}
R_1 &:= \dis\lim_{s\to 1^+} \Bigg\{ \chi(\delta;Q)
\left[\sump\frac{1}{Q^s} - \sump\frac{1}{Q'^s}\right] + \chi(\delta;Q')
\left[\sump \frac{1}{Q'^s} - \sump \frac{1}{Q^s}\right]\Bigg\}
\\[\jot]
&= \dis\lim_{s\to 1^+} \Bigg\{
\left[1 - \biggl( \frac{2}{\delta}\biggr)\right]
\left[\chi(\delta;Q)\sump\frac{1}{Q^s} - 
\chi(\delta;Q')\sump\frac{1}{Q'^s}\right] \Bigg\} \,.
\end{align*}
\end{lemma}

\begin{proof}
We begin with the ``composition identity'':
$$
2AX^2 + CY^2 \equiv (Ax_1^2 + 2Cy_1^2)\, (2x_2^2 + ACy_2^2),
$$
where
$$
X := x_1x_2 + Ay_1y_2,  \qquad  Y := 2x_1x_2 - Cy_1y_2.
$$
This identity \textit{reveals the relation which every form $Q$
has with its homologue}.

We associate the Jacobi symbol
$$
\biggl(\frac{\delta}{2 + AC}\biggr)
$$
to the form
$$
2x_2^2 + ACy_2^2.
$$

But, by the generalized law of quadratic reciprocity for Jacobi
symbols:
$$
\biggl(\frac{\delta}{2 + AC}\biggr)
= (-1)^{\bigl(\frac{\delta-1}{2}\bigr)
\cdot \bigl(\frac{2 + AC - 1}{2}\bigr)}
\biggl(\frac{2 + AC}{\delta}\biggr)
$$
and, since
$$
\delta\equiv 1 \pmod 4  \phrase{and}  \delta \mid AC,
$$
we conclude that
$$
(-1)^{\frac{\delta-1}{2}} = 1
$$
and
$$
2 + AC \equiv 2 \pmod \delta,
$$
which mean
$$
\biggl(\frac{\delta}{2 + AC}\biggr) = \biggl(\frac{2}{\delta}\biggr),
$$
and therefore
$$
\biggl(\frac{\delta}{2A + C}\biggr)
= \biggl(\frac{2}{\delta}\biggr)\biggl(\frac{\delta}{A + 2C}\biggr),
\qquad
\biggl(\frac{\delta}{A + 2C}\biggr)
= \biggl(\frac{2}{\delta}\biggr)\biggl(\frac{\delta}{2A + C}\biggr),
$$
and we obtain
\begin{align}
R_1 &= \dis \lim_{s\to 1^+} \Bigg\{ \left[\chi(\delta;Q) - 
\biggl(\frac{2}{\delta}\biggr)\chi(\delta;Q)\right]
\sump\frac{1}{Q^s} +  \left[\chi(\delta;Q') - 
\biggl(\frac{2}{\delta}\biggr) \chi(\delta;Q') \right]
\sump\frac{1}{Q'^s}     \Bigg\}
\nonumber \\[\jot]
&= \dis \lim_{s\to 1^+} \Bigg\{ 
\left[1 - \biggl(\dis\frac{2}{\delta}\biggr)\right]
\left[\chi(\delta;Q)\sump\frac{1}{Q^s} - 
\chi(\delta;Q')\sump\frac{1}{Q'^s}\right] \Bigg\} \,.
\tag*{\qed}
\end{align}
\hideqed
\end{proof}

Now we compute $R_1$ according as $\biggl(\dfrac{2}{\delta}\biggr)$ is
$+1$ or $-1$, respectively.

\begin{lemma}
$\dis\biggl(\frac{2}{\delta}\biggr) = +1 \implies  R_1 = 0$, and
$$
\biggl(\frac{2}{\delta}\biggr) = -1 \implies
R_1 = 2 \,\lim_{s\to 1^+} \Bigg\{ \chi(\delta;Q)\sump\frac{1}{Q^s}
- \chi(\delta;Q')\sump\frac{1}{Q'^s} \Bigg\}.
\eqno \qed
$$
\end{lemma}

Now we used the value of $R_1$ to compute the \textit{full} weighted
sum~$R$.

\begin{lemma}
$\dis \biggl(\frac{2}{\delta}\biggr) = +1 \implies  R = 0$, and
$$
\biggl(\frac{2}{\delta}\biggr) = -1 \implies
R = \frac{4\pi}{\sqrt{m}} \sum_{k=1}^{h} \chi(\delta,Q_{k})
\ln g_{m/A_k^2} \,.
$$
where the sum is over all $h(-4m)$ reduced forms
$Q_k \equiv A_k X^2 + C_k Y^2$ of determinant~$m$.
\end{lemma}

\noindent
\textbf{Note:}\quad
Case 1 of Proposition 2 explains the first ``miraculous cancellation''
in the computation of~$R$.

\begin{proof}
By Kronecker's formula for $g_n$, we conclude
$$
R = \sum_{k=1}^{h} \chi(\delta,Q_{k})\, \frac{4\pi}{\sqrt{m}}
\ln g_{m/A_k^2}.
$$
But, by Lemma~\ref{lm:R-one},
$$
R = \left[1 - \biggl(\frac{\delta}{2}\biggr)\right] \cdot
\lim_{s\to 1^+} \left\{
\sum_{k=1}^{h}\chi(\delta,Q_{k}) \cdot \sump\frac{1}{Q_{k}^s}
\right\} \,.
$$
If $\dis \biggl(\frac{\delta}{2}\biggr) = 1$, then $R = 0$, as
claimed in the statement.
If $\dis \biggl(\frac{\delta}{2}\biggr) = -1$, then
\begin{align}
R &= 2 \cdot \lim_{s\to 1^+} \left\{
\sum_{k=1}^{h}\chi(\delta,Q_{k}) \cdot \sump\frac{1}{Q_{k}^s}
\right\}
\nonumber \\
&= \frac{4\pi}{\sqrt{m}}\sum_{k=1}^{h}\chi(\delta,Q_{k})
\ln g_{m/A_k^2} \,.
\tag*{\qed}
\end{align}
\hideqed
\end{proof}

\subsubsection{The second cancellation: Dirichlet's computation}

We propose to compute the sum
$$
F := \sum_{\delta} \chi(\delta,Q),
$$
where $\delta$ runs over all $2^{t}$ fundamental discriminants which
divide $D = -4m$ and $Q\equiv AX^2 + BY^2$ is any one of the $2^{t}$
reduced forms of determinant $m$. We will considere two cases: when
$Q$ is the so-called \textit{principal} form
$$
X^2 + mY^2,
$$
and the case when $Q$ is a \textit{non-principal} form and therefore
$$
Q \neq X^2 + mY^2.
$$

\vspace{6pt}
\noindent
\fbox{\textbf{Case 1:  $Q \equiv X^2 + mY^2$}}
\vspace{6pt}

Then
\begin{align*}
\chi(\delta,Q) &= \biggl(\frac{\delta}{m + 1}\biggr)
\\
&= (-1)^{(\frac{\delta-1}{2}) \cdot (\frac{m + 1 - 1}{2})}
\biggl(\frac{m + 1}{\delta}\biggr)
\\
&= \biggl(\frac{1}{\delta}\biggr)
\phrase{since} m + 1 \equiv 1 \pmod \delta,
\\
&= 1.
\end{align*}
Therefore,
$$
F = \sum_{\delta} \chi(\delta,Q) = 2^t = h(m).
$$

\vspace{6pt}
\noindent
\fbox{\textbf{Case 2:  $Q \not\equiv X^2 + mY^2$}}
\vspace{6pt}

Then there exists at least  one fundamental discriminant,
$\delta_1$, for which
$$
\chi(\delta_1,Q) \neq 1
$$
(why?) and thus
\begin{align*}
\chi(\delta_1,Q) \cdot F
&= \chi(\delta_1,Q)\sum_{\delta}\chi(\delta,Q)
\\
&= \sum_{\delta} \chi(\delta_1,Q)\chi(\delta,Q)
\\
&= \sum_{\delta} \chi(\delta_1\,\delta,Q),
\end{align*}
where we have used the multiplicativity of the Jacobi symbol. Now
comes the fundamental observation (no pun intended!): \textit{As
$\delta_1$ runs over the complete set of fundamental discriminants
dividing} $-4m$, \textbf{\emph{so does}} $\delta_1\,\delta$,\emph{ and therefore the
set of Jacobi symbols} $\{\chi(\delta_1,Q)\}$ \textbf{\emph{coincides}} \emph{with the
set of Jacobi symbols} $\{\chi(\delta_1\,\delta,Q)\}$. Therefore,
$$
\chi(\delta_1,Q) \cdot F = F.
$$
But
$$ 
\chi(\delta_1,Q) = -1,
$$
or $-F = F$, i.e., 
$$
F = 0.
$$
Written out,
$$
F := \sum_{\delta} \chi(\delta,Q) = 0
\phrase{for any}  Q \neq X^2 + m Y^2.
$$

\noindent
\textbf{Note:}\quad
Case 2 \textit{explains the second ``miraculous cancellation'' in the
computation of $R$.}

\subsubsection{The general formula for \boldmath $g_m$ where $m$ is even}

If we conjoin the pair of ``miraculous cancellations'' into a single
generic computation we obtain a general formula for Ramanujan's
function~$g_m$.

Our starting point, then, is to compute:
\begin{align*}
H &:= \sum_{\delta} R(\delta)
\\
&= \sum_{\delta} \sum_{k=1}^{h} \chi(\delta,Q_{k}) \frac{4\pi}{\sqrt{m}}
\ln g_{m/A_k^2}
\\
&= \frac{4\pi}{\sqrt{m}} \sum_{k=1}^{h} 
\biggl\{ \sum_{\delta}\chi(\delta,Q_{k}) \biggr\} \ln g_{m/A_k^2}
\\
&=
\begin{cases} 0, &\phrase{if} \biggl(\dfrac{2}{\delta}\biggr) = +1, \\
\dis\frac{4\pi}{\sqrt{m}} \cdot h \cdot \ln g_m, 
&\phrase{if} \biggl(\dfrac{2}{\delta}\biggr) = -1. \end{cases}
\end{align*}
But, as we saw in \eqref{eq:four-eqns},
\begin{align*}
R &= 2 \cdot \lim_{s\to 1^+} \left\{
\sum_{k=1}^{h} \chi(\delta,Q_{k}) \cdot \sump\frac{1}{Q_{k}^s}
\right\}
\\
&= \dis 2 \cdot \left\{\sum_{n=1}^{\infty} 
\biggl( \frac{\delta}{n}\biggr)\frac{1}{n^s} \right\}
\left\{\sum_{m=1}^{\infty}\biggl( \frac{\delta}{m}\biggr)\frac{1}{m^s}
\right\} \,.
\end{align*}

Now we assume that
$$
\delta\delta' = -4m ,
$$
which means that $\delta'$ \textit{is even and of opposite sign to
$\delta$}, and we conclude from the Dirichlet class number formulas
\eqref{eq:Diri-formula-neg} and \eqref{eq:Diri-formula-pos} that
\begin{align*}
\dis \frac{4\pi}{\sqrt{m}} \sum_{k=1}^{h} \chi(\delta,Q_{k})
\ln g_{m/A_k^2}
&= \dis \frac{2\pi}{\sqrt{m}} \cdot K(\delta) \cdot K(\delta')
\cdot \ln \varepsilon
\\
\implies 2 \cdot \sum_{k=1}^{h} \chi(\delta,Q_{k})
\ln g_{m/A_k^2}
&= K(\delta) \cdot K(\delta') \cdot \ln \varepsilon
\\
\implies 2 \cdot \sum_{\delta}\sum_{k=1}^{h} \chi(\delta,Q_{k})
\ln g_{m/A_k^2}
&= \sum_{\delta}K(\delta) \cdot K(\delta') \cdot \ln \varepsilon
\\
\implies 2 \cdot h \cdot \ln g_m
&= \sum_{\delta} K(\delta) \cdot K(\delta') \cdot \ln \varepsilon
\\
\implies (g_m)^{2h}
&= \prod_{\delta} \epsilon^{K(\delta) \cdot K(\delta')} \,,
\end{align*}
or, putting
$$
m = 2n
$$
and
$$
h(m) \equiv K(-4m) = K(-8n),
$$
we obtain our crowning result:
\begin{theorem}Let $\delta_{1}$, $\delta_{2},\cdots,\delta_{K(-8n)}$ the distinct odd divisors of $-8n$, and suppose the complimentary divisor $\delta'_{k}$ to $\delta_{k}$ is defined by\begin{equation*}\fbox{$\dis\delta_{k}\cdot\delta'_{k}:=-8n.$}\end{equation*}  Then

\begin{center}
\Large
\fbox{\fbox{\parbox{10cm}{$\dis
g_{2n} = \prod_{k=1}^{K(-8n)}
\left(\frac{T_{k}+U_{k}\sqrt{\delta_{k}}}{2}\right)^{\frac{K\left(\delta_{k}\right) \cdot K\left(\delta'_{k}\right)}{K(-8n)}},
$} }}
\end{center}
\noindent
\end{theorem}
\noindent a formula which can only be called glorious!

\vspace{6pt}

\noindent
\textbf{Note:}\quad
This formula has been proved by Weber \cite{w1,w2}, and by Berndt
\cite{berndt}. Weber's development is based on quadratic forms, but
invokes \textit{genus} theory, which we wished to avoid. Berndt's
proof uses ideal theory and appeals to results in various journals,
although both he and Weber obtain a more general result. Our
development reveals the fundamental structure of all of these proofs
without the extra baggage of (mathematically useful, but pedagogically
obfuscatory) generalization.

This same formula can be used to compute $g_{2n}$ for all even
convenient numbers which are of the form $2P$ where $P$ is a product
of distinct odd primes. In principle, this is just what Weber did.


\hfill Received \date{\today}

\noindent MARK B. VILLARINO \\
 \indent       Escuela de Matem\'atica,\\
  \indent      Universidad de Costa Rica,\\
  \indent      2060 San Jos\'e, Costa Rica\\
  \indent      mvillari@cariari.ucr.ac.cr\\

\end{document}